\documentclass[11pt]{article}
\usepackage{amssymb}
\usepackage{latexsym}
\usepackage{pb-diagram}
\usepackage{amsmath}
\usepackage{amscd}
\usepackage{oldgerm}
\usepackage{here}

\def\Z{\mathbb{Z}}
\def\R{\mathbb{R}}
\def\C{\mathbb{C}}

\def\fhol{\mathfrak{hol}}
\def\fm{\mathfrak{m}}
\def\fh{\mathfrak{h}}

\def\fso{\mathfrak{so}}

\def\fsu{\mathfrak{su}}
\def\fsp{\mathfrak{sp}}
\def\fg{\mathfrak{g}}
\def\fspin{\mathfrak{spin}}

\def\1{\mathbf{1}}
\def\:{\lrcorner}
\def\#{\sharp}

\def\g{\gamma}

\def\o{\circ}

\def\<#1,#2>{\langle#1,\,#2\rangle}

\def\S{\mathbb{S}\,}

\def\qed{\ensuremath{\quad\Box\quad}}

\def\inv#1{\raise.1em\hbox to 0pt{$^{-1}$\hss}_{#1}\;}

\def\v{\noindent}

\newcommand{\lan}{\langle}
\newcommand{\ran}{\rangle}

\newcommand{\bean}{\begin{eqnarray*}}
\newcommand{\eean}{\end{eqnarray*}}
\newcommand{\benu}{\begin{enumerate}}
\newcommand{\eenu}{\end{enumerate}}
\newcommand{\eea}{\end{eqnarray}}
\newcommand{\bea}{\begin{eqnarray}}

\newcommand{\tvarphi}{\widetilde{\varphi}}
\newcommand{\tpsi}{\widetilde{\psi}}

\newcommand{\tE}{\tilde{E}}

\newcommand{\tW}{\tilde{W}}
\newcommand{\tX}{\tilde{X}}

\newtheorem{Theorem}{Theorem}
\newtheorem{Proposition}{Proposition}

\newtheorem{Corollary}{Corollary}
\newtheorem{Lemma}{Lemma}
\newcommand{\vphi}{\varphi}
\newcommand{\tvphi}{\widetilde{\varphi}}

\title{Codazzi spinors and globally hyperbolic manifolds with special holonomy}

\begin{document}

\author{Helga Baum\footnote{Institut f\"ur Mathematik,
Humboldt-Universit\"at Berlin, Sitz: Rudower Chaussee 25,D-12489
Berlin. email: baum@mathematik.hu-berlin.de } and Olaf
M\"uller\footnote{Instituto de Matem\'aticas,
Universidad Nacional Aut\'onoma de M\'exico (UNAM)
Campus Morelia, C. P. 58190, Morelia, Michoac\'an, Mexico. email: olaf@matmor.unam.mx}}

\date{\today}
\maketitle
\begin{abstract}
\v In this paper we describe the structure of Riemannian manifolds
with a special kind of Codazzi spinors. We use them to construct
globally hyperbolic Lorentzian manifolds with complete Cauchy surfaces, for any weakly irreducible
holonomy representation with parallel spinors, t.m. with a holonomy
group $G \ltimes \R^{n-2} \subset SO(1,n-1)$,
where $G \subset SO(n-2)$ is trivial or a product of groups $SU(k)$, $Sp(l)$, $G_2$ or $Spin(7)$. \\

\end{abstract}

\tableofcontents

\section{Introduction}

The connected holonomy groups of Riemannian manifolds are well
understood and there are a lot of results and methods for
construction of Riemannian metrics with special holonomy (cf.
\cite{jo2}). Contrary to that, the classification of holonomy
groups for indefinite metrics is a longtime and widely open
problem, since the existence of indecomposable but non-irreducible
holonomy representations makes the classification difficult.
Recently, the classification of the connected holonomy groups of
Lorentzian manifolds was achieved. L.~Berard-Bergery and
A.~Ikemakhen described the structure of weakly-irreducible,
non-irreducible subgroups of the Lorentz group (cf. \cite{BI} or
\cite{G1}). Th.~Leistner (\cite{tl1}, \cite{tl2}) was able to
classify all Lorentzian Berger algebras, which is the essential
part in the classification of connected Lorentzian holonomy
groups, and he realized part of them as holonomy algebras.
A.~Galaev (\cite{G2}) gave local analytic Lorentzian metrics for
all of these Berger algebras, including that of the still missing
coupled types, thereby completing the classification of connected
Lorentzian holonomy groups. The next task in this line is the
construction of global geometric models
with special Lorentzian holonomy. \\[0.1cm]
There is a class of manifolds that is very suitable for
purposes of field theories, mathematical physics etc: the class of
globally hyperbolic manifolds, which can be seen as a
Lorentzian analog to complete manifolds in Riemannian
geometry. The construction of analytic manifolds in \cite{G2} of course yields immediately the existence of globally hyperbolic metrics for these holonomies as the holonomy of an analytic manifold equals its local holonomy, and in each analytic Lorentzian manifold every point has a globally hyperbolic neighborhood. But we can sharpen the requirements a bit and try to construct {\em globally hyperbolic manifolds with complete Cauchy surface} which property clearly is not shared by the examples constructed by the method described above. A class of globally hyperbolic manifolds with even stronger
completeness conditions (which imply e.g. strong statements about the
long-time existence of Lorentzian minimal surfaces, cf. \cite{m})
is the one of bbc manifolds which will be defined in section
\ref{sec-cylinders}.
Therefore we ask for bbc manifolds with special holonomy.
In this article we give a
construction method which yields such manifolds using
an idea of Ch.~B\"ar, P.~Gauduchon and A.~Moroianu (\cite{bgm}),
who constructed parallel spinors on generalized cylinders out of
Codazzi spinors. We explore this method for the Lorentzian
situation in detail, describe the structure of all Riemannian
manifolds with imaginary Codazzi spinors to an invertible Codazzi
tensor as well as the causal and holonomy properties of the
Lorentzian cylinder defined by the Codazzi tensor
(which will be a bbc manifold).\\[0.1cm]
In section \ref{sec-cylinders} we start with basic properties of
Lorentzian cylinders by which we mean a Lorentzian manifold $(a,b)
\times M$ with a metric $g$ of the form $ g=- dt^2 + g_t$, where
$g_t$ is a smooth family of Riemannian metrics on $M$. In §2.1 we
give a criterion for global hyperbolicity of a Lorentzian cylinder.
In §2.2 we describe basics of Lorentzian spin geometry in order to
fix our conventions and notations for the rest of the paper. In
§2.3 we consider special Lorentzian cylinders that are constructed
out of Codazzi spinors on the Riemannian base $(M,g_0)$. A spinor
field $\vphi$ on $(M,g_0)$ is called {\em imaginary Codazzi
spinor} if
\[ \nabla^{g_0}_X \vphi = i A(X) \,\small{\bullet}\, \vphi \]
for all vector fields $X$, where $A$ is a Codazzi tensor on
$(M,g_0)$, which is then uniquely determined by the spinor
$\varphi$. As in \cite{bgm} we prove, that $\vphi$ induces a
parallel spinor $\tvphi$ on the Lorentzian cylinder
$(C:=(a,b)\times M\,,\; g:=-dt^2 + (1-2tA)^*g_0)$ and analyze the
causal type of the associated Dirac current on $C$. In section
\ref{sec-invCodazzi} we draw our attention to the case of
invertible Codazzi tensors. In § 3.1 we prove that an imaginary
Codazzi spinor with invertible Codazzi tensor $A$ to the metric
$g_0$ corresponds to an imaginary Killing spinor for the metric
$A^*g_0$. Applying the classification result for manifolds with
imaginary Killing spinors we obtain the structure result for
manifolds admitting imaginary Codazzi spinors with invertible
Codazzi tensor $A$. Any such complete manifold is isometric to
\[ M= \R \times F\,,\; g_0= (A^{-1})^*(ds^2 + e^{-4s}g_F) \qquad \qquad (*) \]
where $(F,g_F)$ is a complete Riemannian manifold with parallel
spinors and $A^{-1}$ is a Codazzi tensor on the warped product $\R
\times_{e^{-2s}} F$ (Theorem \ref{Th-Codazzi-spinor1}). In § 3.2
 we analyze the existence of
Codazzi tensors on warped products and reduce it to the question
of existence of Codazzi tensors on $(F,g_F)$. This leads to a
construction method for Riemannian manifolds with imaginary
Codazzi spinors. Any $(n-1)$-dimensional Riemannian manifold
$(F,g_F)$ with parallel spinors and Codazzi tensor $T$ which has
eigenvalues uniormally bounded from below gives rise to a Codazzi
tensor $H$ on the warped product $(M=\R \times F, g_{wp}=ds^2 +
e^{-4s}g_F)$ with eigenvalues bounded away from zero. Hence on
$(M,g_0=H^*g_{wp})$ there are Codazzi spinors and on the cylinder
$(C,g_C)$ with
\bean C&=& C[F;H]= (a,b) \times M = (a,b)\times \R \times F \\
g_C&=& -dt^2 + g_t = -dt^2 + (H-2t\1)^*g_{wp}\eean parallel
spinors by the construction explained above.
 In § 3.3 we
study the causal type and the holonomy of this cylinder. In
Theorem \ref{th-cylinder1} we prove that $C[F;H]$ is globally
hyperbolic if $(F,g_F)$ is complete.
Furthermore, we show that
$C[F;H]$ is flat if and only if $(F,g_F)$ is flat. Note, that the
holonomy group of a simply connected Lorentzian manifold
$\widetilde{C}$ acts irreducible if and only if it is isomorphic
to $SO_0(T_x\widetilde{C})$. Since we have a parallel spinor on
the cylinder $C[F;H]$ which defines a nontrivial parallel vector
field by its Dirac current, $C[F;H]$ can not be irreducible. Hence
the cylinder $C[F;H]$ is either weakly irreducible, meaning that
the holonomy representation has no {\em non-degenerate} invariant
subspace (but a degenerate one) or $C[F;H]$ is decomposable  t.m.
it is (locally) the product of a Lorentzian and a Riemannian
manifold. In fact, $C[F;H]$ is a Brinkman space, t.m. it admits a
nowhere vanishing parallel light-like vector field. In Theorem
\ref{Theorem-hol1} we prove that $C[F;H]$ is decomposable if
$(F,g_F)$ contains a flat factor. If $(F,g_F)$ is (locally) the
Riemannian product of irreducible (non-flat) manifolds, then
$C[F;H]$ is weakly irreducible and the connected component of the
holonomy group of $C[F;H]$ is isomorphic to \[ Hol^0(F,g_F)
\ltimes \R^{n-1} . \] This gives a construction method for
globally hyperbolic (and even bbc) manifolds with complete Cauch hypersurfaces for every weakly-irreducible,
non-irreducible Lorentzian holonomy representation with a fixed spinor.
In section 4 we finish the paper by studying
examples of Riemannian manifolds with parallel spinors and Codazzi
tensors which are the building blocks of our construction. \\

\section{Lorentzian cylinders}\label{sec-cylinders}

\subsection{Definition and causal properties}

{\bf Definition:} A {\em Lorentzian cylinder} is a product
manifold $C=(a,b) \times M$ with a Lorentzian metric of the form
$g = -dt^2 \oplus g_t$ where $g_t$ is a smooth family of
Riemannian metrics on $M$ parametrized over the interval $(a,b)$,
where $-\infty \leq a < 0 < b \leq + \infty$.\\

We want to link this notion to the notion of bbc manifolds as
described in \cite{m}:\\[0.2cm]
{\bf Definition:} A Lorentzian manifold $(C,g)$ is said to be {\em
bbc} iff
\[(C,g) \cong (\R \times M, g(p) = -f(p) d t^2 \oplus
g_{t}), \] where $f$ is a smooth positive function on $C$ which is
{\em bounded} on every $\{ t \} \times M$, and $g_t$ is a family
of {\em complete} Riemannian metrics on $M$, and with the
additional property that the eigenvalues of $\dot{g}_t \o g_t^{-1}
$ are uniformally {\em bounded} for all $t$.\\[0.2cm]
\v {\em Examples}: If $g_t $ is a compact perturbation
(e.g. described by pullback by a symmetric endomorphism
which is the identity outside of a compactum on every
$\{ t \} \times M$) of a fixed complete Riemannian
metric for all values of $t$, then for any bounded $f$,
$(\R \times M, - f dt^2 \oplus g_t)$ is bbc. On the contrary,
any open proper subset of the Minkowski space is not bbc as
it does not contain any complete spacelike hypersurface.

\begin{Proposition} (\cite{m})
The $t=constant$ hypersurfaces of a bbc manifold are Cauchy hypersurfaces.
In particular, bbc manifolds are globally hyperbolic.
\end{Proposition}
\bigskip

\v We start with the causality properties of a cylinder. For that,
let $A_t$ denote the connecting endomorphism between $g_t$ and
$g_0$:

\[g_t = A_t^*g_0 = g_0(A_t \,\cdot\, , \,A_t \,\cdot\,). \]

\begin{Proposition}
\label{Prop-geocylinder} Let $(C,g) \cong ( (a, b)  \times M, g =
- d t^2 \oplus g_{t}) $ be a Lorentzian cylinder. Then $(C,g)$ is
strongly causal. If in addition the metric $g_0$ is complete and
the eigenvalues of $(A_t)^{-1}$ are uniformally bounded on $M$ for
every $t \in (a,b)$, then $(C,g)$ is globally hyperbolic. If in
addition to both previous assumptions the eigenvalues of
$\dot{A}_t \o A_t^{-1}$ are uniformally bounded on every $\{ t \}
\times M$, then $(C,g)$ is bbc.
\end{Proposition}

\v {\bf Proof.} The function $f(t,x)=t$ is obviously a time
function on (C,g), since it is strictly increasing along any
future directed causal curve. Therefore, the Lorentzian cylinder
is stably causal, which implies strong causality (cf. \cite{bee},
chapt. 3).

\v Now assume additionally that the metric $g_0$ is complete and
that for every $t \in (a,b)$ there exists a constant $c_t \in
\R^+$ such that all eigenvalues of $A_t^{-1}$ are uniformally
bounded by $c_t<+\infty$. We want to prove global hyperbolicity of
$C$. Following \cite{bee}, p.65, knowing that we have strong
causality for $C$, we have to show that the intersection
$N:=J^+(y_0) \cap J^-(y_1)$ of the causal future $J^+(y_0)$ and
the causal past $J^-(y_1)$ is compact for all points $y_0 =
(T_0,q_0)$ and $y_1 = (T_1, q_1)$ of $C$. For that, let $\{x_n\}$
be a sequence of points in $N$. Then $x_n=(t_n,p_n)$, where $t_n
\in [T_0,T_1]$. Hence $\{t_n\}$ has a convergent subsequence. So,
we may assume that $\{t_n\}$ converges to $t^* \in [T_0,T_1]$.
Now, let $\;\gamma_n(s)=(t(s),\delta_n(s)): [0,1] \rightarrow N
\;$ be a future directed causal curve in $N$ from $y_0$ to $x_n$.
Then by the assumption on $A_t$
\[ c \|\delta_n'(s)\|_{g_0} < \|\delta_n'(s)\|_{g_{t(s)}} \leq t'(s)
\]
for a constant $c > 0$. It follows that the length of $\delta_n$
with respect to the metric $g_0$ is bounded by
$R=\frac{1}{c}(T_1-T_0)$. Hence, all points $p_n$ lie in the
$g_0$-geodesic ball $B^{g_0}(q_0,R)$ of radius $R$ around $q_0$.
Since $g_0$ is complete, this ball is relatively compact, so
$\{p_n\}$ has a convergent subsequence. This proves that $N$ is
compact.

\bigskip

\v The bbc property follows directly from the definition as $f =1$ and $g_t = A_t ^* g_0 $ is complete by the assumption on the eigenvalues of $A_t^{-1}$ before. This completes the proof.   \hfill \qed


\subsection{Spinors on Lorentzian cylinders}

In this section we describe spinors on a Lorentzian cylinder. For
convenience of the reader we first recall some basic facts about
Lorentzian spin geometry, thereby fixing our notations.
For details we refer to \cite{baum1} and \cite{bgm}.\\
Let $(\R^{1,n},\eta)$ be the $(n+1)$-dimensional Minkowski space
with the inner product
\[ \eta(x,y)= -x_0y_0 + x_1y_1 + x_2y_2 + \dots + x_ny_n,\]
where $x=(x_0,x_1,\dots,x_n)$, $y=(y_0,y_1,\dots,y_n)\,$. We fix
the standard isometric embedding of the Euclidian space
$(\R^n,\epsilon)$ into the Minkowski space
\[
\begin{array}{cccc} i: &(\R^n,\epsilon)&  \longrightarrow
& (\R^{1,n},\eta)\\
& x & \longmapsto & (0,x)
\end{array}
\]
with the timelike normal vector $e_0=(1,0,\dots,0)$. \\
We denote by $Cl_n:=Cliff(\R^n,\epsilon)$ and
$Cl_{1,n}:=Cliff(\R^{1,n},\eta)$ the Clifford algebras of the
Euclidian space and the Minkowski space, respectively. Since the
complex linear map
\[ \begin{array}{cccc}\beta: & (\C^n,\epsilon^{\C}) & \longrightarrow &
Cliff^0(\C^{1,n},\eta^{\C}) \\
& x & \longmapsto & i e_0\cdot x
\end{array}
\]
satisfies
\[ \beta(x) \cdot \beta (x) = - \epsilon^{\C}(x,x) \cdot 1 , \]
$\beta$ induces an isomorphism of the complexified Clifford
algebras
\[ \tau: Cl_n^{\C} \longrightarrow (Cl_{1,n}^0)^{\C}. \]
Now, let us consider the usual Spin-representation of $Cl_{1,n}$
on the spinor modul $\Delta_{1,n}$
\[ \rho_{1,n}: Cl_{1,n} \longrightarrow GL(\Delta_{1,n}) \]
and let denote by
\[\rho_{1,n}^{\pm}: Cl_{1,n} \longrightarrow GL(\Delta_{1,n}^{\pm}) \]
its positive and negative parts in case of odd $n$.\\
If $n$ is even, we consider the action of the Clifford algebra
$Cl_n$ on the space $\Delta_{1,n}$ given by
\[ \begin{array}{ccccccl} \kappa:= \rho_{1,n}  \circ \tau: & Cl_n &
\longrightarrow &
(Cl_{1,n}^0)^{\C} & \longrightarrow & GL(\Delta_{1,n})& \\
& x & \longmapsto &  i e_0\cdot x & \longmapsto & \rho_{1,n}(i \
e_0 \cdot x)& ,\qquad x \in \R^n
\end{array}
\]
which is a realization of the unique irreducible complex
representation of $Cl_n$. Let us denote by $\Delta_n$ the
$Spin(n)$-representation \bea \Delta_n: =
(\Delta_{1,n},\kappa|_{Spin(n)}).  \label{spinors1} \eea Then as
$Spin(n)$-representations $\Delta_n$ and
$(\Delta_{1,n},\rho_{1,n}|_{Spin(n)})$ are isomorphic and the
Clifford product under this isomorphism is identified by
\bea x \star u \in \Delta_n \longmapsto i \ e_0 \bullet x \bullet
u \in \Delta_{1,n} \; \qquad x \in \R^n, \;u \in \Delta_n
,\label{spinors2} \eea where $x\star u := \kappa(x)u$ denotes the
Clifford product on $\Delta_n$ and $x \bullet u:= \rho_{1,n}(x)u$
denotes the Clifford product on $\Delta_{1,n}$.\\
If $n$ is odd, the actions of $Cl_n$ on
$\Delta_{1,n}^{\pm}$ given by
\[ \begin{array}{ccccccl} \kappa^{\pm}:= \rho_{1,n}^{\pm}  \circ
\tau: & Cl_n & \longrightarrow &
(Cl_{1,n}^0)^{\C} & \longrightarrow & GL(\Delta_{1,n}^{\pm})& \\
& x & \longmapsto &  i e_0\cdot x & \longmapsto & \rho_{1,n}^{\pm
}(i\ e_0 \cdot x) & ,\qquad x \in \R^n
\end{array}
\]
realize the two non-equivalent irreducible complex representations
of $Cl_n$.\\
We denote by $\Delta_n$ and $\hat{\Delta}_n$ the
$Spin(n)$-representations \bea \Delta_n &:=&
(\Delta_{1,n}^+,\kappa^+|_{Spin(n)}) \nonumber \\
\hat{\Delta}_n & := & (\Delta_{1,n}^-,\kappa^-|_{Spin(n)})
 \label{spinors1-a}
\eea The $Spin(n)$-representations $(\Delta_n,\kappa^+)$ and
$(\Delta_{1,n}^+,\rho_{1,n}^+)$ and the
$Spin(n)$-repre\-sen\-ta\-tions $(\hat{\Delta}_n,\kappa^-)$ and
$(\Delta_{1,n}^-,\rho_{1,n}^-)$ are isomorphic, whereby the
Clifford product is identified by
\bea x \star u &:=& \kappa_+(x)u \in \Delta_n \longmapsto i \ e_0
\bullet x \bullet u \in \Delta_{1,n}^+ \;
\qquad x \in \R^n, \;u \in \Delta_n ,\label{spinors2-a}\\
 x \star v &:=& \kappa_-(x)v \in \hat{\Delta}_n \longmapsto i \ e_0 \bullet x \bullet
 v
\in \Delta_{1,n}^- \; \qquad x \in \R^n, \;v \in \hat{\Delta}_n.
\label{spinors2-b} \eea Furthermore, the linear isomorphism
\[ \begin{array}{cccc} \phi: &  \Delta_n = \Delta_{1,n}^+ &
\longrightarrow & \hat{\Delta}_n = \Delta_{1,n}^-  \label{spinors2-c} \\
& u & \longmapsto & \hat{u}:= e_0 \bullet u
\end{array}
\]
is an isomorphism of the $Spin(n)$-representations and the
Clifford product on $\Delta_n$ and $\hat{\Delta}_n$ satisfies \bea
\widehat{x \star u} = - x \star \hat{u}\;, \qquad x \in \R^n, \; u
\in \Delta_n=\Delta_{1,n}^+ .\label{spinors2-d}
 \eea

On the Lorentzian spinor modul $\Delta_{1,n}$ we have two
hermitian products $\lan \cdot,\cdot \ran_1$ and $\lan \cdot,\cdot
\ran_0$ which are related by
\bea \lan u,v \ran_1 = \lan e_0 \bullet u,v \ran_0 \;, \qquad u,v
\in \Delta_{1,n} . \label{spinors3} \eea
 The inner product $ \lan \cdot,\cdot \ran_1$ is
$Spin_0(1,n)$-invariant and the Clifford product $\,x \bullet\;$
is symmetric for all vectors $x \in \R^{1,n}$. The inner product
$\lan \cdot,\cdot \ran_0$ is $Spin(n)$-invariant and the Clifford
product
$\,x\star\; $ is skew-symmetric for all vectors $x \in \R^n$.\\[0.1cm]
Now, let $(M,g_0)$ be a $n$-dimensional Riemannian manifold and
let $(C=(a,b)\times M,g=-dt^2+g_t)$ be a Lorentzian cylinder over
$(M,g_0)$. In the following we denote by $\nu$ the  timelike unit
normal field on $C$ given by
\[ \nu(t,x):= \frac{\partial}{\partial t}(x,t) . \]
Now, let us assume that $(M,g_0)$ is a spin manifold with the
$SO(n)$-frame bundle $P_M$, the spin structure $(Q_M,f_M)$ and the
corresponding spinor bundles
\bean S_M & :=&  Q_M \times_{Spin(n)} \Delta_n \\
\hat{S}_M & :=& Q_M \times_{Spin(n)} \hat{\Delta}_n \; \qquad
\mbox{if $n$ is odd }.\eean

The spin structure $(Q_M,f_M)$ induces a spin structure on the
Lorentzian cylinder $(C,g)$ in a canonical way. To explain that,
consider $\gamma_x(t):=(t,x)$ which is a timelike geodesic through $(0,x)$
and let us denote by
\[ \tau_{\gamma_x}^t : T_{(0,x)}(\{0\}\times M) \longrightarrow
T_{(t,x)}(\{t\}\times M) \] the parallel displacement in $(C,g)$
along $\gamma_x$ with respect to the Levi-Civita connection of $g$
and by $\pi: C \longrightarrow M$ the projection $\pi(t,x)=x$.
Then the $SO(n)$-principal bundle
\[\begin{array}{ll} \hat{P}:= \Big\{ (\nu(t,x), \tau^t_{\gamma_x}(s_1),\dots,
\tau^t_{\gamma_x}(s_n)) \mid  & x \in M, \; (s_1,\dots,s_n) \; \mbox{positively}\\
& \mbox{oriented ON-basis in $(T_xM,(g_0)_x)$} \Big\}
\end{array} \]
is a $SO(n)$-reduction of the $SO_0(1,n)$-principal bundle $P_C$ of oriented and time-oriented frames on
$(C,g)$ with $\pi^*P_M \simeq \hat{P}$. Therefore the $SO_0(1,n)$-frame bundle
of $(C,g)$ can be described as
\[
P_C= \pi^*(P_M \times_{SO(n)} SO_0(1,n)) \] and the following pair
$(Q_C,f_C)$ is a spin structure of $(C,g)$
\[Q_C := \pi^* (Q_M \times_{Spin(n)} Spin_0(1,n)) \]
\[\begin{array}{rccc}
 f_C: &  Q_C & \longrightarrow & P_C \\
 &[q,A] & \longmapsto &  [f_M(q),\lambda(A)] \end{array}
 \]
where $\, \lambda: Spin_0(1,n) \longrightarrow SO_0(1,n)\,$
denotes the
usual 2-fold covering. \\
Using (\ref{spinors1}) and (\ref{spinors1-a}) we obtain the
following identification for the spinor bundles $S_C$ of $(C,g)$
and $S_M$ of $(M,g_0)$. For even $n$ hold
\bea S_C &:=& Q_C \times _{(Spin_0(1,n),\rho)}\Delta_{1,n} \simeq
\pi^*(Q_M \times_{(Spin(n),\kappa)} \Delta_n) \nonumber\\
& = & \pi^* S_M \label{spinors-cylinder1} \eea and for odd $n$
 \bea S_C &=& S_C^+ \oplus S_C^- \nonumber\\
& =&  Q_C \times _{(Spin_0(1,n),\rho^+)} \Delta_{1,n}^+ \; \oplus
\; Q_C \times _{(Spin_0(1,n),\rho^-)} \Delta_{1,n}^- \nonumber
\\
& \simeq & \pi^*\Big(Q_M \times_{(Spin(n),\kappa^+ )} \Delta_n
\,\oplus\,
Q_M \times_{(Spin(n),\kappa^-)}\hat{\Delta}_n\Big) \nonumber \\
& =&  \pi^* S_M \oplus \pi^* \hat{S}_M \label{spinors-cylinder1-a}
\eea

Hence, for even $n$, any spinor field $\psi \in \Gamma(S_C)$ on
the cylinder can be understood as an $1$-parameter family of
spinors $t \in (a,b) \mapsto \psi(t,\cdot) \in \Gamma(S_M)$ on the
manifold $M$. In case of odd $n$, any spinor field in
$\Gamma(S_C^+)$ is a $t$-parameter family of spinor fields in
$\Gamma(S_M)$, whereas any spinor field in $\Gamma(S_C^-)$ can be
understood as $t$-parameter family of spinors in $\Gamma(\hat{S}_M)$.    \\
Now, for any vector field $X$ and any spinor field $\varphi$ on
$(M,g_0)$ we denote by $\tX$ and $\tvarphi$ the vector field and
the spinor field on $(C,g)$ arising from $X$ and $\varphi$ by
parallel displacement along the geodesics $\gamma_x$
\[ \tX(t,x):= \tau_{\gamma_x}^t(X(x)) \quad ,\quad \tvarphi(t,x):=
\tau_{\gamma_x}^t(\varphi(x)) .\]
Because of (\ref{spinors2}) and (\ref{spinors2-a}) and  the
Clifford multiplication satisfies \bea \widetilde{X \star \varphi}
= i \, \nu \bullet \tX \bullet \tvarphi
\label{spinors-cylinder2}.\eea
Due to its invariance properties, the hermitian products $\lan
\cdot ,\cdot \ran_1$ and $ \lan \cdot,\cdot \ran_0$ on the spinor
modul $\Delta_{1,n}$ induce hermitian inner products on the spinor
bundel $S_C$ which are related by
\bea \lan \psi_1,\psi_2 \ran_1 = \lan \nu \bullet \psi_1,\psi_2
\ran_0 \;, \qquad \psi_1,\psi_2 \in \Gamma(S_C).
\label{spinors-cylinder3} \eea
For any spinor field $\psi \in \Gamma(S_C)$ and $\varphi \in
\Gamma(S_M)$ we define the associated Dirac current $V_{\psi}$ on
$(C,g)$ and $W_{\varphi}$ on $(M,g_0)$, respectively, by
\bea g(V_{\psi},Y) &:=& - \lan Y \bullet \psi, \psi \ran_1 \qquad
\forall Y \in {\cal X}(C)  \label
{spinors-cylinder4}\\
g_0(W_{\varphi},X) &:=& \; i \lan X \star \varphi, \varphi \ran_0 \qquad
\forall X \in {\cal X}(M) \label{spinors-cylinder5}\eea The Dirac
currents satisfy the following conditions, which are easily to
verify using (\ref{spinors2}), (\ref{spinors2-a}),
(\ref{spinors-cylinder4}) and (\ref{spinors-cylinder5})

\begin{Proposition}\label{prop-spinors1}
\begin{enumerate}
\item The Dirac current $V_{\psi}$ of a nowhere vanishing spinor field
 $\psi$ on the Lorentzian
cylinder $(C,g)$ is a causal and future directed vector field.
\item Let $\varphi$ be a spinor field on $(M,g_0)$ and $\tvarphi$
its parallel extension to the cylinder $(C,g)$ along the geodesics
$\gamma_x$. Then for the Dirac currents hold \bea V_{\tvarphi} =
\|\tvarphi\|_0^2 \;\nu + \widetilde{ W_{\varphi}}
\label{spinors-cylinder6} \eea
\end{enumerate} \hfill \qed
\end{Proposition}

\v Comparing the Levi-Civita connections of the cylinder $(C,g)$
and of the level sets $(M_t:=\{t\}\times M ,g_t)$ one obtains the
following relation between the spinor derivatives:

\begin{Proposition}\label{prop-spinors2}(cf. \cite{bgm})
Let $X$ be a vector field and $\varphi$ a spinor field on
$(M,g_0)$. Then for the parallel transported spinor field
$\tvarphi$ on $(C,g)$ hold $\nabla^C_{\nu}\tvarphi = 0 $ by definition and
\bea \nabla^{C}_{\tX}\tvarphi &=&
\nabla^{M_t}_{\tX}\tvarphi + \frac{1}{2} \nu \bullet S_t(\tX))
\bullet \tvarphi  \qquad \mbox{on $M_t$}\label{spinors-cylinder7}\eea

\v where $S_t(Y):= - \nabla^C_Y \nu $ is the Weingarten map of the
submanifold $M_t \subset C$. \hfill \qed
\end{Proposition}
\ \\


\subsection{Codazzi spinors and special Lorentzian cylinders}

If $B$ is an invertible endomorphism field and $h$ a
pseudo-Riemannian metric on $M$,  we denote by $B^*h$ the induced
metric\[ B^*h(X,Y):= h(BX,BY). \] In the following we will
consider a special example of
Lorentzian cylinders.\\
Let $A$ be a symmetric uniformally bounded endomorphism field. We
denote by $\mu_+(A)$ the supremum of the positive eigenvalues of
$A$ or zero if all eigenvalues of $A$ are nonpositive and by
$\mu_-(A)$ the infimum of the negative eigenvalues of $A$ or zero
if all eigenvalues of $A$ are nonnegative. Now, set
 $a:=
(2\mu_-(A))^{-1}$ and $b:=(2\mu_+(A))^{-1}$. We denote by $C(M;A)$
the Lorentzian cylinder
\[ C(M;A) := (a,b) \times M \;, \quad g:= -dt^2 + g_t := -dt^2 + (\1-2tA)^*g_0\]
call it {\em
Lorentzian cylinder over $M$ induced by $A$}.\\

{\bf Definition:} Let $(M,g_0)$ be a Riemannian spin manifold with
spinor bundle $S_M$. An {\em imaginary Codazzi spinor} on $M$ is a
spinor field $ \psi \in \Gamma(S_M)$ which satisfies the equation
\bea \nabla_X \psi = i \cdot A(X) \star \psi , \label{Codazzi1}
\eea
where $A$ is a Codazzi tensor on $(M,g_0)$, t.m. $A$ is a symmetric
endomorphism field which obeys the Codazzi equation
$$(\nabla^M _X A) (Y) = (\nabla^M_Y A) (X) $$
for all vector fields $X,Y$ on $M$. The equation (\ref{Codazzi1}) can be interpreted as the one of parallelity w.r.t. the (non-metric) connection $\nabla -i A$. \\

The notion of imaginary Codazzi spinors generalizes the case of imaginary Killing spinors with Killing number $i\lambda$, which are Codazzi spinors with
$A=\lambda \1$.\\
Note that the Codazzi-tensor $A$ is uniquely defined by its
Codazzi spinor (\ref{Codazzi1}) since
\[ g_0(A(X),Y) = -\frac{1}{2}\; Im\Big(\;\lan X \star \nabla_Y^M \psi + Y
\star \nabla^M_X \psi , \psi \ran_0\Big) \cdot \|\psi\|_0^{-2}.\]
Furthermore, the existence of an imaginary Codazzi spinor implies
the following curvature constraint on $(M,g_0)$

\begin{Proposition}\label{Prop-Codazzi-spinors1}
Let $(M,g_0)$ be a Riemannian manifold with non-trivial imaginary
Codazzi spinor $\psi$ to the Codazzi tensor $A$. Then the Ricci
curvature of $(M,g_0)$ satisfies
\[ Ric^M(X) = 4 A^2(X) - 4 (trA)\cdot A(X) .\]
\end{Proposition}

\v {\bf Proof.} The Codazzi equation for $A$ implies
\[ R^{S_M}(X,Y)\psi = \Big( A(X) \star A(Y) - A(Y)\star A(X)\Big) \star
\psi .\] Hence\\[-0.7cm]
 \bean  Ric(X) \star \psi=  -2 \; \sum_{k=1}^{n} s_k
\star R^{S_M}(X,s_k) \psi  = (4 A^2(X) - 4 tr(A)\cdot A(X)) \star
\psi. \eean Since $\psi$ is nontrivial, it vanishes nowhere.
Therefore, the vectors in front of $\psi$ in the latter formula on
both sides coincide.  \hfill \qed
\ \\

 In \cite{bgm} the authors
study {\em real} Codazzi spinors and cylinders with {\em
spacelike} cylinder axis and show, that Codazzi spinors on $M$ can
be extended to parallel spinors on the cylinder. The same
statement is true in case of imaginary Codazzi spinors and
Lorentzian cylinders with the completely analogous proof:

\begin{Proposition}\label{Prop-codazzi-spinors2}
Let $(M,g_0)$ be a Riemannian spin manifold carrying an
imaginary Codazzi spinor $\psi \in \Gamma(S_M)$ with uniformally
bounded Codazzi tensor $A$ and let $C:=C(M;A)$ be the Lorentzian
cylinder over $(M,g_0)$ induced by $A$ with its canonical spin
structure. Then the $\nu$-parallel extension $\tpsi \in
\Gamma(S^{(+)}_{C})$ of $\psi$ is a parallel spinor field on the
cylinder $C(M;A)$. Conversely, the restriction of any parallel
spinor $\phi \in \Gamma(S^{(+)}_{C})$ of the cylinder to $M_0$ is
a Codazzi spinor to the Codazzi tensor $A$.
\end{Proposition}

\v {\bf Proof.} We only recall the main steps of the proof and refer
for details to \cite{bgm}. \\
First, it is easy to see that $2A$ is the Weingarten map of the
submanifold $M_0 = \{0\} \times M \subset C(M;A)=:C$. Then from
(\ref{spinors-cylinder2}) and (\ref{spinors-cylinder7}) follows
that for any vector field $X$ on $M_0$ \bea \nabla^C_X \tpsi =
\nabla^{M_0}_X \tpsi + \nu \bullet A(X)\bullet \tpsi =
\nabla^{M}_X \psi - i A(X) \star \psi . \label{Cod}\eea
 Since $\psi$ is an
imaginary Codazzi spinor, we obtain $\; \nabla^C_X \tpsi = 0\;$
for all vectors $X$ tangential to $M_0$. The Codazzi
 equation for $A$ implies that $\nu$ is in the kernel of the
 curvature endomorphism, i.e. $R^C (\nu , \cdot) = 0$. Since
 $[\nu,X]=0$ and $\nabla_{\nu}^C \tpsi =0 $ it follows
 \[\nabla^{C}_{\nu}
  \nabla^{C}_{X} \tpsi = \frac{1}{2} R^C (\nu, X) \bullet \tpsi = 0 \]
where $X$ is the lift of a vector field $X$ of $M$ to $C$. But
this means that the spinor field $\nabla^C_{X}\tpsi $ is parallel
along all geodesics $\gamma_x$. As this spinor field vanishes in
the point $\gamma_x(0)=(0,x)$, it vanishes everywhere on $C$.
Hence $\tpsi$ is parallel on $C$. On the other hand, any parallel
spinor $\phi$ on the cylinder is the parallel extension of its
restriction to $M_0$ along the curves $\gamma_x$. (\ref{Cod})
shows that this restriction is an imaginary Codazzi tensor. \hfill
\qed
\ \\

For a spinor field $\varphi \in \Gamma(S_M)$ we consider the
 the real subbundle $\; E_{\vphi}:= TM \star \vphi \subset S_M$
 and denote by $dist_{\vphi}$ the pointwise distance between $i \vphi$ and
 $E_{\vphi}$ with respect to the real inner product $\;Re \lan \cdot
 , \cdot \ran_0\;$ on $S_M$. \\
In order to analyze the causal type of the Dirac current
$V_{\tpsi}$ of the parallel extension $\tpsi$ of a Codazzi spinor
$\psi$ we will use the following Lemma

\begin{Lemma}\label{Lemma-Codazzi1}
Let $(M,g_0)$ be a Riemannian manifold with an endomorphism field
$B$ and let $\varphi\in \Gamma(S_M)$ be a spinor field such that
\bea \nabla^M_X \varphi = i\, B(X) \star \varphi.
\label{Codazzi2}\eea We denote by $W_{\varphi}$ the Dirac current
of $\vphi$ and by $q_{\vphi}$ the function \bea q_{\varphi} :=
\|\vphi\|^4_0 - g_0(W_{\vphi},W_{\vphi}). \label{Codazzi3} \eea
Then $q_{\vphi}$ is constant and non-negative and given by the
distance function
 \bean q_{\vphi} = dist_{\vphi}^2 \cdot \|\vphi\|_0^2.
 \eean
\end{Lemma}

\v {\bf Proof.} Using formula (\ref{spinors-cylinder5}) for the Dirac
current $W_{\vphi}$ and formula (\ref{Codazzi2}) for the spinor
$\vphi$ we obtain
\[ \nabla^M_X W_{\vphi} = 2 \|\vphi\|_0^2 \cdot B(X). \]
From that follows
 \bean X(q_{\vphi}) &=& 4 \|\vphi\|_0^2 \;i\; \lan
B(X)\star \vphi,\vphi \ran_0
- 2 g_0(\nabla_X W_{\vphi},W_{\vphi})\\
&=& 4 \|\vphi\|_0^2 g_0(B(X),W_{\vphi}) - 4 \|\vphi\|_0^2
g_0(B(X),W_{\vphi})
 = 0.\eean
If $(s_1,\dots,s_n)$ is an ON-basis on $M$, then $\;
(\|\vphi\|^{-1}\,s_j \star \vphi \mid j=1,\dots,n)\;$ is an
ON-basis in the real vector space $E_{\vphi}$. Hence
\bean
dist_{\vphi}^2 &=& \|i\vphi\|_0^2 - \|proj_{E_{\vphi}}(i\vphi)\|_0^2 \\
&=& \|\vphi\|_0^2 - \|\sum_{j=1}^{n} \|\vphi\|_0^{-2} \cdot \lan
i\vphi , s_j \star \vphi \ran_0 \cdot s_j \star \vphi \|_0^2\\
&=& \|\vphi\|^2_0 - g_0(W_{\vphi},W_{\vphi}) \cdot
\|\vphi\|_0^{-2} \;=\; q_{\vphi}\cdot \|\vphi\|_0^{-2}.
 \eean
\hfill \qed

\bigskip

\v Now, we describe the causal type of the Dirac current of a
parallel spinor on the cylinder which is induced by an imaginary
Codazzi spinor.

\begin{Proposition}\label{Prop-codazzi-spinors3}
Let $\psi$ be a non-vanishing imaginary Codazzi spinor on
$(M,g_0)$ with uniformally bounded Codazzi tensor $A$ and let
$\tpsi$ be its $\nu$-parallel extension to the cylinder $C(M;A)$.
If $\;dist_{\psi}=0\;$, the Dirac current $V_{\tpsi}$ of $\tpsi$
is parallel and lightlike. If $\;dist_{\psi}
> 0\;$, the Dirac current $V_{\tpsi}$ is parallel and timelike.
\end{Proposition}
{\bf Proof.} Since $\tpsi$ is parallel, the Dirac current
$V_{\tpsi}$ is parallel as well. $\psi$ is non-trivial, hence the
length $\|\psi\|_0^2$ has no zeros. From formula
(\ref{spinors-cylinder6}) follows \bean
g(V_{\tpsi},V_{\tpsi}) &=& -\|\tpsi\|^4_0 +
g(\tW_{\psi},\tW_{\psi}) = -\|\psi\|^4 + g_0(W_{\psi},W_{\psi}) = -
q_{\psi} \\ &=& - \|\psi\|^2_0 \cdot dist_{\psi}
 \eean
 which proves the statement. \hfill \qed
 \ \\

We remark, that for any spinor $\vphi$ on a $3$- and on a
$5$-dimensional manifold the distance $dist_{\vphi}$ is zero (cf. \cite{BFGK}, p. 89). \\


\section{Codazzi spinors with invertible Codazzi tensors}\label{sec-invCodazzi}

\subsection{The structure theorem}

 Let us denote by $d^{\nabla}$ the exterior derivative
induced by a covariant derivative $\nabla$. In particular, for any
endomorphism field $A$
\[(d^{\nabla} A)(X,Y) = (\nabla_X A)(Y) - (\nabla_Y A)(X). \]
Hence $A$ is a Codazzi tensor with respect to $\nabla$ if and only
if $d^{\nabla}A = 0$.\\
\v In the following we want to examine the case of invertible Codazzi
tensors more closely. First we remark

\begin{Proposition}
\label{Prop-connections} (\cite{dc}) Let $A$ be an invertible
Codazzi tensor on $(M,g_0)$ and let $\nabla^{g_0}$ and
$\nabla^{A^*g_0}$ denote the Levi-Civita connections of $g_0$ and
of the induced metric $A^*g_0$. Then
\[\nabla^{A^*g_0}_X  = A^{-1} \circ \nabla^{g_0}_X \circ A \;,
\qquad \forall \; X \in {\cal X}(M) \]

\v and $A^{-1}$ is a Codazzi tensor w.r.t. $A^*g_0$. In
particular, the holonomy groups of $g_0$ and $A^*g_0$ are
conjugated
\[ Hol_x(M,A^*g_0) = A_x^{-1} \circ Hol_x(M,g_0) \circ A_x \]
and for the curvature hold
\bean R^{A^*g_0}(X,Y,Z,V) &=& R^{g_0}(X,Y,AZ,AV) \\
Ric^{g_0} \circ A &=& A \circ Ric^{g_0}.  \eean \hfill \qed
\end{Proposition}

 \v There is an obvious identification between the spin
structures of the manifolds $(M,g_0)$ and $(M,A^*g_0)$. Let us
denote by
\[ \Phi_A : \; \varphi \in S_M= S_{(M,g_0)} \longrightarrow \overline{\varphi} \in
\overline{S}_M:=S_{(M,A^*g_0)}\]

\v the corresponding isomorphism between the spinor bundles. Then
for the Clifford product and the spinor derivatives
$\nabla=\nabla^{g_0}$ and $\overline{\nabla}= \nabla^{A^*g_0}$
hold \bea \overline{X \star \vphi} = A^{-1}(X) \star
\overline{\vphi}\qquad \mbox{and}\qquad \overline{\nabla_X \vphi}
= \overline{\nabla}_X \overline{\vphi}.\label{bar1}\eea
Furthermore, we obtain \bea \lan \vphi,\psi \ran_0 = \lan
\overline{\vphi},\overline{\psi} \ran_0 \qquad \mbox{and} \qquad
W_{\overline{\vphi}} = A^{-1}(W_{\vphi}) \label{bar2}\eea which
shows that $q_{\vphi} = q_{\overline{\vphi}}$ for the functions
$q_{\vphi}$ defined in (\ref{Codazzi3}).
 Hence we obtain the following Proposition:

\begin{Proposition}\label{Prop-Killingspinors1}
Let $\vphi$ be an imaginary Codazzi spinor on a Riemannian
manifold $(M,g_0)$ with invertible Codazzi tensor $A$. Then the
corresponding spinor $\overline{\vphi}$ on $(M, A^*g_0)$ is a
Killing spinor with Killing number i. Vice versa, given a
Riemannian manifold $(M,g_0)$ with a Killing spinor $\phi$ to the
Killing number $i a$ for some real number $a$ and with an
invertible Codazzi tensor $A$, on the Riemannian manifold $(M,
A^*g_0)$ the spinor field $\Phi_A(\phi)=\overline{\phi}$ is a
Codazzi spinor with Codazzi tensor $aA^{-1}$. Likewise, $\Phi_A$
 maps parallel spinor fields (if they exist) to parallel spinor
 fields.
 \hfill \qed
\end{Proposition}

\v Thus we can use the structure results the first author proved
for manifolds with imaginary Killing spinors.

\begin{Proposition}(\cite{h1}, \cite{h2})\label{Prop-Killingspinors2}
\v  $(M,g_0)$ is a complete connected Riemannian spin manifold
with non-vanishing imaginary Killing spinor to the Killing number
$i\mu$ if and only if it is isometric to a warped product
\[ (\R \times F, ds^2 + e^{-4\mu s}g_F) \]
where $(F,g_F)$ is a complete connected spin manifold with
non-vanishing parallel spinor.  \hfill \qed
\end{Proposition}

\v Let $B$ a symmetric, positive definite endomorphism field on $(M,g_0)$ whose eigenvalues are
uniformally bounded away from zero, t.m. with $0 < \|B^{-1} \| \leq c < \infty $ on all of $M$.
Then the metric $B^*g_0$ is complete if $g_0$ is complete. From
Proposition \ref{Prop-Killingspinors1}, Proposition
\ref{Prop-Killingspinors2} and Proposition
\ref{Prop-codazzi-spinors3} we obtain the following structure of a
Riemannian manifold with imaginary Codazzi spinors.

\begin{Theorem}\label{Th-Codazzi-spinor1}
Let $(M,g_0)$ be a complete connected Riemannian spin manifold
with non-vanishing imaginary Codazzi spinor to a Codazzi tensor
$A$ whose eigenvalues are uniformally bounded away from zero. Then
$(M,A^*g_0)$ is isometric to a warped product
\[ (\R \times F, \; g_{wp}:= ds^2 + e^{-4 s}g_F) \qquad (*)\]
where $(F,g_F)$ is a complete connected spin manifold with
non-vanishing parallel spinor. Furthermore, $A^{-1}$ is a Codazzi
tensor on the warped product {\em (*)}.  \hfill \qed
\end{Theorem}
\ \\[-0.5cm]
\subsection{Codazzi tensors on warped products}

Theorem \ref{Th-Codazzi-spinor1} shows how we can obtain
Riemannian manifolds with Codazzi spinors. For that we have to
find Codazzi tensors on warped products of the form
\[ M = \R \times_f F := ( \R \times F, g_M = ds^2 + f^2(s) g_F).  \qquad \qquad \]

Let us now address to this question. We decompose $T_{(s,x)}M
\simeq \R
\partial_s \oplus T_x F$ and write a symmetric $(1,1)$-tensor
field $H$ on $M$ with respect to this decomposition as
$$ H = H (b,D,E) =
\left(
\begin{array}{cc}
b\cdot Id &  \widetilde{D}\\
D & E
\end{array} \right)$$

where $b$ is a real function on $\R \times F$, $E$ is a
$s$-parametrized family of symmetric endomorphism fields on
$(F,g_F)$, $D$ an endomorphism field from $\R \partial_s$ to $TF$,
and

\[
 \widetilde{D}V
= f^2\cdot g_F(V,D(\partial_s)) \,\partial_s\; \]
 for all vector fields
$V$ on $F$.\\

\v Let us call the tensor field $H=H(b,0,E)$ {\em simple} iff
$E(s) = K(s) \cdot Id_F$ for all $s \in \R$.

\begin{Proposition}\label{Prop-Codazzi4}
The endomorphism field $H(b,D,E)$ is a Codazzi tensor on the
warped product $\;M= \R \times_f F\;$ if and only if for all
vector fields $V$ and $W$ on $F$ hold
 \bea
\nabla^F_V(D(\partial_s)) &=& \dot{E}(V) + \frac{\dot{f}}{f} E(V)
- b \frac{\dot{f}}{f}V \label{Codazzi5}\\
 grad^Fb &=& 3 f \dot{f} D(\partial_s) + f^2 \dot{D}(\partial_s)
 \label{Codazzi6}\\
 (d^{\nabla^F}E)(V,W) &=& f \dot{f}\Big( g_F(V,D(\partial_s))W -
 g_F(W,D(\partial_s))V \Big) \label{Codazzi7}
 \eea In particular, the
$s$-parameter family $D(\partial_s)$ of vector fields on $F$
satisfies \bea R^F(V,W)D(\partial_s) = (f \ddot{f} - \dot{f}^2)
\Big(g_F(V,D(\partial_s))W - g_F(W,D(\partial_s))V \Big).
\label{Codazzi8}\eea Here $\dot{}$ denotes the derivative with
respect to the parameter $s$.
\end{Proposition}

{\bf Proof.}  We rewrite $d^{\nabla^M}H = 0$ in conditions for
$b$, $D$ and $E$. The covariant derivative on the warped product
$M= \R \times_f F$ is given by
\bean
\nabla^M_{\partial_s} \partial_s &=& 0 \\
\nabla^M_{\partial_s}V &=& \nabla^M_V \partial_s \,=\,
\frac{\dot{f}}{f} V \\
\nabla^M_V W &=& - f \dot{f} g_F(V,W) \partial_s + \nabla^F_V W.
\eean (Cf. \cite{ONeill}, p. 206)). Using these formulas, one
easily sees that the formulas (\ref{Codazzi5}) and
(\ref{Codazzi6}) are equivalent to $\;
(d^{\nabla^M}H)(\partial_s,V) = 0\;$. The condition
$\;(d^{\nabla^M}H)(V,W)= 0\;$ is equivalent to (\ref{Codazzi7})
and
\[ g_F(W,\nabla^F_V D(\partial_s)) = g_F(V, \nabla^F_W
D(\partial_s)),\] where the latter formula follows already from
(\ref{Codazzi5}). Now, if we differentiate (\ref{Codazzi5}) again
with respect to $\nabla^F_W$ and insert (\ref{Codazzi6}) and
(\ref{Codazzi7}) we obtain (\ref{Codazzi8}). \hfill \qed

\ \\

Recall that a Codazzi tensor is called {\em trivial} if it is a
constant multiple of the Identity.

\begin{Corollary}\label{Cor-Codazzi1}

\begin{enumerate}
\item On every warped product $M=\R
\times_f F$ there is a nontrivial invertible Codazzi tensor.
Explicitly, if $T$ is a (possibly trivial) Codazzi tensor on
$(F,g_F)$, then we obtain a Codazzi tensor on $M$ by $\,H =
H(b,0,E)\,$ where $b$ is a function which depends only on $s$ and
$E$ is given by \bea E (s) = \frac{1}{f} \Big(T + \int_0^s
b(\sigma) \dot{f}( \sigma ) d \sigma \cdot Id_{F}
\Big)\label{Codazzi9}\eea
\item Let $f(s)= e^{k \cdot s}$. If there is a non-simple Codazzi tensor on
$M$ then there is a nontrivial Codazzi tensor on $F$ or a nonzero, parallel
or homothetic vector field on $F$.
\end{enumerate}
\end{Corollary}

\v {\bf Proof.} It is easy to check, using (\ref{Codazzi5}),
(\ref{Codazzi6}) and (\ref{Codazzi7}), that the tensor field
$H(b,0,E)$, where $b=b(s)$ is a real function and $E$ is given by
(\ref{Codazzi9}), is a Codazzi tensor on $M=\R \times_f F$, for
every Codazzi tensor $T$ on $F$. If we choose for example $b$ so that
$b\dot{f}>0$ and $T=k\cdot Id_F$ with $k>0$, then $H(b,0,T)$ is invertible. \\
 Now let $H(b,D,E)$ be a
Codazzi tensor on $M=\R \times_f F$, first with a general warping
function $f$. Let us first consider the case that $D\equiv 0$.
Then formula (\ref{Codazzi6}) shows, that the function $b$ depends
only on $s$. From (\ref{Codazzi7}) follows, that $E$ is a Codazzi
tensor on $(F,g_F)$. From (\ref{Codazzi5}) we obtain the following
differential equation in the space of Codazzi tensors on $(F,g_F)$
\[ \dot{(fE)} = \dot{f} E + f \dot{E} =  b \dot{f} Id_F. \]
Hence $E$ is given by formula (\ref{Codazzi9}), where $T$ is a
Codazzi tensor on $(F,g_F)$. $H$ is simple exactly if $T$ is
trivial.\\
Now, let $D\not = 0$. We set
\[\tE := (f\ddot{f}-\dot{f}^2) E - f \dot{f}  \nabla^F
D(\partial_s).\] From (\ref{Codazzi8}) follows
\[
(d^{\nabla^F}d^{\nabla^F}D(\partial_s))\,(V,W)
=R^F(V,W)D(\partial_s) = (f \ddot{f} - \dot{f}^2)
\Big(g_F(V,D(\partial_s))W - g_F(W,D(\partial_s))V \Big).
\]
Using (\ref{Codazzi7}) we obtain
 \bean
(d^{\nabla^F}\tE) (V,W)&=& (f\ddot{f}-\dot{f}^2) (d^{\nabla^F}E)
(V,W) - f \dot{f}\;
(d^{\nabla^F}d^{\nabla^F}D(\partial_s))\,(V,W) \\
&=& (f\ddot{f}-\dot{f}^2) \Big( (d^{\nabla^F}E)(V,W) - f \dot{f}
(g_F(V,D(\partial_s))W-g_F(W,D(\partial_s))V ) \Big)  \\
&=& 0. \eean Thus $\tE$ is a Codazzi tensor on $(F,g_F)$.

In our case, for $f=e^{k \cdot s}$, we have $f\ddot{f}-\dot{f}^2
\equiv 0$. Therefore, $\tilde{E} = f \dot{f} \nabla^F D
(\partial_s) $, thus $\nabla^F D (\partial_s) $ is a Codazzi
tensor on $(F,g_F)$. If this Codazzi tensor is trivial, t.m.
 $\nabla^F_{\cdot} D(\partial_s) =  \lambda Id_F$, then
$D(\partial_s)$ is a family of parallel vector fields on $F$ (if
$\lambda=0$) or of homothetic vector fields on $F$ (if $\lambda
\not = 0$). \hfill \qed

\ \\[-0.5cm]
\begin{Corollary}\label{Cor-Codazzi2}
Let $M=\R \times_f F$ be a warped product with warping function
$f(s)=e^{-2s}$ and let $T$ be a Codazzi tensor on $(F,g_F)$ whose eigenvalues are uniformally bounded from below. Then the
warped product $M$ admits a Codazzi tensor $H$ with eigenvalues
uniformally bounded away from zero ($ \|H ^{-1}\| \leq c < \infty$ on all of $M$).
\end{Corollary}

{\bf Proof.} Let all eigenvalues of $T$ be greater then $k\in \R$.
We choose a strictly increasing function $h\in C^{\infty}(\R)$
with $h(s) < \frac{k}{2}$. Let $b(s):= e^{2s}\dot{h}(s)$. Then all
eigenvalues of the Codazzi tensor $H=H(b,0,E)$, where $E$ is given
by (\ref{Codazzi9}), are nonnegative. Hence the norm of the
Codazzi tensor $H + c Id_M $ for a positive constant $c$ is
bounded from below by $c$. \hfill \qed

\ \\
\subsection{The holonomy and causal properties of the special Lorentzian cylinder}

Theorem \ref{Th-Codazzi-spinor1} and Corollary \ref{Cor-Codazzi2}
provide a construction principle for Riemannian manifolds with
Codazzi spinors and for Lorentzian manifolds with special
holonomy. We start with a connected Riemannian manifold $(F,g_F)$
which admits a non-vanishing parallel spinor and a Codazzi tensor
with eigenvalues uniformally bounded from below. Then the warped
product $\,(M:=\R \times F, \,g_{wp}:= ds^2 + e^{-4s}g_F)\,$
admits imaginary Killing spinors to the Killing number $i$ and a
Codazzi tensor $H$ with positive eigenvalues uniformally bounded away from zero.
The Riemannian manifold $(M,g_0:=H^*g_{wp})$ admits imaginary
Codazzi spinors for the Codazzi tensor $H^{-1}$.

In the following we will denote the Lorentzian cylinder
constructed in this way out of $(F,g_F)$ and $H$ by $C[F;H]$. By
construction
\[ C[F;H]:= (a,b)\times M = (a,b) \times \R \times F\]
where $a$ is the half of the supremum of all negative eigenvalues
of $H$ or $-\infty$ if all eigenvalues of $H$ are positive and $b$
is the half of the infimum of all positive eigenvalues of $H$ or
$+\infty$ if all eigenvalues of $H$ are negative. The metric of
$C[F;H]$ is given by
\[ g_C := -dt^2 + g_t = -dt^2 + (H-2t\1)^*(ds^2 + e^{-4s}g_F). \]

In the following we will denote
\bean A &:=& H^{-1}\\
 H_t &:= &(H-2tId_M)
 \eean

First, let us start with the causality property of the cylinder
$C[F;H]$.

\begin{Theorem}\label{th-cylinder1}
If the Riemannian manifold $(F,g_F)$ is complete, then the
cylinder $C[F;H]$ is globally hyperbolic and even bbc.
\end{Theorem}
{\bf Proof.} Since $(F,g_F)$ is complete, the warped product
$(M=\R\times_{e^{-2s}} F, g_{wp})$ is complete as well. The eigenvalues
of $(H-2t\1)$ are uniformally bounded away
from zero for all $t\in (a,b)$, hence Proposition
\ref{Prop-geocylinder} yields that the Lorentzian cylinder
$C[F;H]$ is globally hyperbolic and moreover bbc \hfill   \qed

\ \\ Next we would like to calculate the curvature of $C[F;H]$.
For that we note

\begin{Lemma}\label{lem-P1}
The Weingarten map $W_t$ of the hypersurface $M_t=\{t\} \times M
\subset C[F;H]$ is given by \bea W_t(X) = 2 H_t^{-1}(X)
\label{P1}.\eea The covariant derivative of $C[F;H]$ is \bea
\nabla^C_X Y &=& ( H_t^{-1}\circ \nabla^{wp}_X\circ H_t) Y -
2g_{wp}(H_tX,Y)
\partial_t \label{P2}\\
\nabla^C_X \partial_t &=& -2 H_t^{-1}X \label{P3}\\
\nabla^C_{\partial_t}X &=& -2H_t^{-1}X \label{P4}\\
\nabla^C_{ \partial_t}\partial_t &=& 0 \label{P5}\eea where $X$
and $Y$ are lifts of vector fields of $M$.
\end{Lemma}
{\bf Proof.} For the Weingarten map $W_t(X) =
-\nabla^C_X\partial_t$ of $M_t$ we have
\[ g_t(W_t(X),Y) = - \frac{1}{2}\dot{g}_t(X,Y).\]
since  $g_t = (H-2tId)^*g_{wp}$ we obtain
\[ \dot{g}_t(X,Y) = -4 g_{wp}(H_tX,Y) = -4 g_t(X,H_t^{-1}Y) =
-4 g_t(H_t^{-1}X,Y) \] Hence
\[ W_t(X) = 2 H_t^{-1}(X). \]
The Gau{\ss} decomposition of the covariant derivative gives
\bean \nabla^C_X Y &=& \nabla^{g_t}_X
Y - g_t(\nabla^C_X Y,\partial_t)
\partial_t = \nabla^{g_t}_X Y -g_t(W_t(X),Y) \partial_t \\
&=& \nabla^{g_t}_X Y -2g_t(H_t^{-1}(X),Y)\partial_t =
\nabla^{g_t}_X Y - 2g_{wp}(X,H_tY)\partial_t \\
&=&\nabla^{g_t}_X Y - 2g_{wp}(H_tX,Y)\partial_t. \eean
Then (\ref{P2}) follows from Proposition \ref{Prop-connections}, since $H_t$ is a Codazzi tensor. The equation (\ref{P5}) holds for any
cylinder metric and (\ref{P3}) and (\ref{P4}) follow from
(\ref{P1}).\hfill \qed

\ \\

\begin{Proposition}
\label{Prop-curvatureFC} The curvature of the Lorentzian cylinder
$C[F;H]$ satisfies:\\[-0.2cm]

1. $R^C(X,Y)\partial_t = R^C(X, \partial_t)Y =
R^C(X,\partial_t)\partial_t = 0$ for all
vector fields $X$, $Y$ on $C$.\\

2. If $X, Y, V \in TM_t$ then
 \bea R^{C}(X,Y)V &=&
(H_t^{-1} \circ R^{wp}(X,Y)\circ H_t)\,V  \nonumber\\
&& - 4 g_{wp}(X,H_tV)H_t^{-1}Y + 4 g_{wp}(Y,H_tV) H_t^{-1}X
\label{curv1}\eea

2. If the vectors $X,Y,V \in TM_t$ are lifts of vectors in $TF$
then  \bea R^{C}(X,Y) H_t^{-1}V = H_t^{-1}  R^F(X,Y)V
\label{curv2} \eea

3. If the vectors $X,Y \in TM_t$ that are lifts of vectors in $TF$
then
 \bea R^{C}(X,Y)H_t^{-1}(\partial_s) = R^{C}(\partial_s,Y)H_t^{-1}(\partial_s)
  = R^C(\partial_s,Y)H_t^{-1}(X) =
 0. \label{curv3}
 \eea
  In particular, $C[F;H]$ is flat if and only if $(F,g_F)$ is flat.
\end{Proposition}

{\bf Proof.} The first statement is a direct calculation using the
Codazzi-Mainardi equation and the Ricatti equation (comp.
\cite{bgm}). The Codazzi-Mainardi equation (\cite{ONeill}, p.115)
shows that $R^C(X,Y)V$ is tangent to $M_t$. Using the Gau{\ss}
equation (\cite{ONeill}, p. 100) and the formula (\ref{P1}) for
the Weingarten map we obtain
\begin{align*} & g_t(R^C(X,Y)V,W) = \\
& \quad =  g_t (R^{M_t} (X,Y) V, W) - 4 ( g_t(H_t^{-1}X,V) g_t
(H_t^{-1}Y,W) + 4 g_t (H_t^{-1}X,W) g_t (H_t^{-1}Y,V))
\end{align*}
Since $\; g_t=H_t^*g_{wp}\;$ we can apply Proposition
\ref{Prop-connections} which yields
\bean
 g_t(R^C(X,Y)V,W) &=&  H_t^* g_{wp} (R^{H_t^*g_{wp}} (X,Y) V, W)
- 4 g_{wp} (X, H_tV) g_t (H_t^{-}Y, W)\\
& & + 4 g_t (H_t^{-1}X, W) g_{wp}(Y, H_tV) \\
& =&  g_t (H_t^{-1}\,R^{g_{wp}} (X,Y) H_tV, W) - 4 g_{wp} (X,
H_tV) g_t(H_t^{-1}Y, W) \\
& & + 4  g_{wp} (Y, H_tV) g_t (H_t^{-1}X, W). \eean
 Hence \bean
R^{C}(X,Y)V = H_t^{-1} \circ R^{wp}(X,Y)\circ H_t \,V - 4
g_{wp}(X,H_tV)H_t^{-1}Y + 4 g_{wp}(Y,H_tV) H_t^{-1}X. \eean \v For
a warped product $M = \R \times_f F$ and vector fields $U,X,Y$
which are lifts of vector fields on $F$, we have
\bea R^{wp} (X,Y) U &=&  R^F (X,Y) U + 4 g_{wp} (X,U) Y - 4
g_{wp} (Y,U) X  \label{P12}\\
R^{wp}(\partial_s,Y)U &=& - 4 f^2 g_F(Y,U) \partial_s \label{P13}\\
R^{wp}(\partial_s,Y) \partial_s &=& 4 Y \label{P14}\\
R^{wp}(X,Y)\partial_s &=& 0. \label{P15}\eea
 \v (cf. \cite{ONeill}, p. 210).
This shows formulas (\ref{curv2}) and (\ref{curv3}). \hfill \qed

\ \\
Finally we would like to study the holonomy of the cylinder
$C[F;H]$. First, let us make a comment on the holonomy of the
manifolds $(M=\R\times F\,,\, H_t^*g_{wp})\,$.

\begin{Proposition}\label{Prop-holM}
The connected component of the holonomy group of the manifolds
$(M=\R\times F\,,\,H_t^*g_{wp})$ is isomorphic to $SO(n)$.
\end{Proposition}
{\bf Proof.} By Proposition \ref{Prop-connections} it is enough to
prove this for $(M,g_{wp})$. Fix a point $z=(0,x) \in M$ and
consider the skew symmetric endomorphism $B:= 4 \partial_s \wedge
Y \subset \fso(T_zM,g_{wp})\simeq \Lambda^2(T_zM)$. Then formulas
(\ref{P13}) and (\ref{P14}) show that \bean
 B(U) &=& -4 g_F(Y,U)\partial_s \;= \; R^{wp}_z(\partial_s,Y)U \\
 B(\partial_s)&=& 4Y \;=\;R^{wp}_z(\partial_s,Y)\partial_s \eean
 for all $Y,U \in T_xF$. By the Ambrose-Singer Theorem
 $\; \partial_s \wedge Y = \frac{1}{4}B = \frac{1}{4} R_z^{wp}(\partial_s,Y) \in \fhol_z(M,g_{wp})\;$
 for all $Y \in T_xF$. Since $[\partial_s \wedge X, \partial_s \wedge
 Y]= X \wedge Y\; $ if follows that
 $\fhol_z(M,g_{wp})=\fso(T_zM,g_{wp})$. \hfill \qed

\ \\[0.2cm]
We will determine the holonomy group of $C[F;H]$ by calculating
the parallel displacement explicitly. For that we use the
following Lemmata.

\begin{Lemma}\label{lem-current}
The Dirac current of any parallel spinor $\tvphi$ on $C[F;H]$ that
is induced by a parallel spinor on $(F,g_F)$ is light-like and
given by
\[ V_{\tvphi}= e^{-2s}(\partial_t - \widetilde{A\partial_s)}\]
where $\widetilde{X}$ denotes the parallel displacement of a
vector field $X$ on $M_0$ along the $t$-lines of the cylinder.
\end{Lemma}
{\bf Proof.} In \cite{h1} and \cite{h2} it is proven, that for any
imaginary Killing spinor $\psi$ on $(M,g_{wp})$, constructed out
of a parallel spinor of $(F,g_F)$, hold $q_{\psi}=0$. The length
and the Dirac current $W_{\psi}$ of $\psi$ are given by
\[ \|\psi(s,x)\|^2_0 = e^{-2s} \quad \mbox{and} \quad W_{\psi}= - e^{-2s}\partial_s \]
where we normalize $\psi$ by $1=\|\psi(0,\cdot)\|_0^2$.  Let us
denote by $\vphi=\overline{\psi}$ the corresponding Codazzi spinor
on $(M,g_0:=H^*g_{wp})$. Then (\ref{bar2}) shows that the length
and the Dirac current of $\vphi$ are given by
\[\|\vphi\|^2_0 = e^{-2s} \quad \mbox{and} \quad W_{\vphi}= -
e^{-2s}H^{-1}(\partial_s).\] According to Proposition
\ref{prop-spinors1} and Proposition \ref{Prop-codazzi-spinors3}
the Dirac current of the parallel extended spinor $\tvphi$ on the
cylinder $C[F,H]$ induced by $\vphi$ is parallel and lightlike and
given by
\[ V_{\tvphi} = e^{-2s}(\partial_t -
\widetilde{H^{-1}(\partial_s)}). \]\hfill \qed
 \ \\

\begin{Lemma}\label{lem-P2}
Let $Z$ be a vector in $TM_0$. Then the parallel displacement
$\widetilde{AZ}$ of $AZ$ along the $t$-lines of the cylinder is
given by \bea\widetilde{AZ} = H_t^{-1}(Z) \label{P6}\eea
If $V$ and $W$ be lifts of vector fields of $F$, then \bea
\nabla^C_{\partial_s}\widetilde{A\partial_s} &=& -2 \partial_t \label{P7}\\
\nabla^C_{\partial_s}\widetilde{AV} &=& -2 H_t^{-1}(V) \label{P8}\\
\nabla^C_{V}\widetilde{A\partial_s} &=& -2 H^{-1}_t(V) \label{P9}\\
\nabla^C_{V}\widetilde{AW} &=& -2e^{-4s}g_F(V,W) \Big(\partial_t -
H_t^{-1}(\partial_s)\Big) + H_t^{-1}(\nabla^F_V W)\label{P10} \eea
\end{Lemma}
{\bf Proof.} From Lemma \ref{lem-P1} we obtain
\[ \nabla^C_{\partial_t}(H_t^{-1}Z) = \nabla^C_{H_t^{-1}Z}
\partial_t + [\partial_t,H_t^{-1}Z] = -2 H_t^{-2}Z + (H_t^{-1})'Z.
\]
Now, from $H_t \circ H_t^{-1}= Id$ follows \[0= H_t' \circ
H_t^{-1} + H_t \circ (H_t^{-1})' = -2 H_t^{-1} + H_t \circ
(H_t^{-1})'\] and therefore
\[ (H_t^{-1})' = 2 H_t^{-2}.\]
Hence $H_t^{-1}Z$ is parallel along the $t$-lines. \\
To prove the formulas for the covariant derivative we use the
formulas (cf. \cite{ONeill}) for the covariant derivative of warped products
\bean \nabla^{wp}_{\partial_s}\partial_s &=& 0 \\
\nabla^{wp}_{\partial_s}V &=& \nabla^{wp}_V \partial_s = -2V \\
\nabla^{wp}_{V}W &=& 2 e^{-4s}g_F(V,W)\partial_s + \nabla^F_V W
\eean and the formulas (\ref{P2}) and (\ref{P6}). \hfill \qed

\ \\[0.3cm]
Now let us denote by $P$ and $Q$ the following light-like
vector fields on $C[F;H]$ \bean P &=& e^{-2s}(\partial_t
-\widetilde{A\partial_s}) \\
Q &=& \frac{1}{2}e^{2s}(\partial_t + \widetilde{A\partial_s})
\eean Then $g_C(P;Q)=1$ and $P$ is parallel (which is clear, since
$P$ is by Lemma \ref{lem-P2} the Dirac current of a parallel
spinor, or can be calculated with the formulas (\ref{P7}) and
(\ref{P9})). The tangent space $TC$ decomposes into the following
subspaces
\[ TC = \R P \oplus \widetilde{ATF} \oplus \R Q. \]
\ \\
\begin{Proposition}\label{Prop-parallel}
Let $\delta(r)=(t(r),s(r),\gamma(r))$ be a curve in $C=C[F;H]$
starting in $\delta(0)=(0,0,x)\in C$. Then the parallel
displacement of $AZ\in AT_xF$ is given by
\[ \tau_{\delta(r)} (AZ) = U(r)= a_{\gamma,Z}(r) P + e^{2s(r)}\widetilde{AY(r)} \]
where $Y(r)= \tau^F_{\gamma|_{[0,r]}}(Z)$ is the parallel
displacement of $Z$ along $\gamma$ in $F$, and $a_{\gamma,Z}$ is
the function determined by
\[ \dot{a}_{\gamma,Z}(r)= 2g_F(Y(r),\dot{\gamma}(r)) \quad ,\quad
a_{\gamma,Z}(0)=0. \]
\end{Proposition}
{\bf Proof.} We consider the vector field
\[ U(r)= a(r) P + b(r)Q + \widetilde{AZ(r)}.\]
Using (\ref{P7}), (\ref{P8}) and (\ref{P9}) we easily check that \bean \nabla^C_{\partial_t} Q & = & 0\\
\nabla^C_{\partial_s} Q &= & 0 \\
\nabla^C_{\dot{\gamma}}Q &=& -2 e^{2s} H_t^{-1}(\dot{\gamma})\eean
Together with (\ref{P10}) this yields
\[ \nabla^C_{\dot{\delta}} U = \dot{a}P + \dot{b}Q - 2be^{2s}H_t^{-1}(\dot{\gamma})
-2\dot{s}H_t^{-1}(Z(r))-2e^{-2s}g_F(\dot{\gamma},Z(r))P +
H_t^{-1}(\nabla^F_{\dot{\gamma}}Z(r)) \] Since $P$, $Q$ and
$H_t^{-1}TF$ are independent, it follows that $U$ is parallel
along $\delta$ if and only if \bean \dot{b}(r)&=& 0\\
\dot{a}(r)&=& 2 e^{-2s}g_F(\dot{\gamma},Z(r))\\
\nabla^F_{\dot{\gamma}}Z(r)& =& 2be^{2s}\dot{\gamma} + 2\dot{s}
Z(r) \eean If we choose the initial conditions $a(0)=b(0)=0$ and
set $Y(r)= e^{-2s}Z(r)$ we obtain, that $U$ is parallel if and
only if
\bean b(r) &=& 0\\
\dot{a} (r)&=& 2 g_F(\dot{\gamma},Y(r))\\
\nabla^F_{\dot{\gamma}}Y(r) &=& 0 \eean This proves the
Proposition. \hfill \qed

\ \\

\begin{Theorem}\label{Theorem-hol1}
Let $Hol_{\hat{x}}(C,g_C)$ be the holonomy group of the cylinder
$C=C[F;H]$ with respect to the point $\hat{x}\in C$, where
$\hat{x}=(0,0,x)$.
 \begin{enumerate}
 \item If $(F,g_F)$ contains a
flat factor $(F_0,g_{F_0})$, then $C[F;H]$ is decomposable.
\item If $(F,g_F)$ splits (locally) into a Riemannian product of
irreducible non-flat factors, then $C[F;H]$ is weakly irreducible
and the connected component of its holonomy group is given by
\[ Hol^0_{\hat{x}}(C,g_C) = \left( H^{-1} \circ Hol^0_x(F,g_F)) \circ
H \right) \; \ltimes \,\R^{dim(F)}.\]
\end{enumerate}
\end{Theorem}

{\bf Proof.} Let us fix an ON-basis $(v_1,\dots,v_{n-1})$ in
$T_xF$ with respect to $g_F$. Now, choose a closed curve
$\delta(r)=(t(r),s(r),\gamma(r))$  in $C$ with
$\delta(0)=\delta(1)=(0,0,x)$. Then for the parallel displacements
hold according to Proposition \ref{Prop-parallel}
\[ \tau^C_{\delta}(AZ) = a_{\gamma,Z}(1)P + A(\tau^F_{\gamma}(Z)) =
\tau_{\gamma}^C(AZ)\] where $a_{\gamma,Z}$ is determined by
$\dot{a}_{\gamma,Z}(r) =
2g_F(\tau^F_{\gamma|_{[0,r]}}Z,\dot{\gamma}(r))$ and
$a_{\gamma,Z}(0)=0$. In particular, with respect to the basis
$\;(P,Av_1,\dots,Av_{n-1},Q)\;$ of $\;T_{\hat{x}}C\;$ we have \bea
\tau^C_{\delta}= \tau^C_{\gamma}= \left( \begin{array}{ccc}
  1 & a_{\gamma} & * \\[0.1cm]
  0 &  A\circ \tau^F_{\gamma} \circ A^{-1} & *  \label{P11}\\[0.1cm]
  0 & 0 & 1
\end{array}\right)
\eea
 where $a_{\gamma}=(a_{\gamma,v_1},
\dots,a_{\gamma,v_{n-1}}))$. (Note, that the two * in the matrix
are uniquely determined by the other entries since
$\tau^C_{\delta}$ is in $SO(T_{\hat{x}}C,g_C)$). This implies that
\[ Hol_{\hat{x}}(C,g_C) \subset \Big(H^{-1}\circ Hol_x(F,g_F) \circ
H\Big)\, \ltimes \R^{dim(F)}.\]
Now, let us suppose, that  $(F,g_F)=
(F_0,g_{F_0}) \times (F_1,g_{F_1})$ is a Riemannian product of a
flat factor $(F_0,g_{F_0})$ and another factor $(F_1,g_{F_1})$.
 Since $(F_0,g_{F_0})$ is flat, we have for the parallel
displacement along a curve $\gamma=(\gamma_0,\gamma_1)$ with
$\gamma(0) = \gamma(1)=x$
\[ \tau_{\gamma}^F = \left(\begin{array}{cc}
  \tau^{F_0}_{\gamma_0} & 0 \\
  0& \tau^{F_1}_{\gamma_1}
\end{array}\right)
= \left( \begin{array}{cc}
  Id & 0 \\
  0 & \tau^{F_1}_{\gamma_1}
\end{array}\right)
\]
Furthermore, let $(v_1=\partial_{x_1},\dots,v_k=\partial_{x_k})$
be a parallel coordinate basis of an Euclidean chart in $F_0$ and
$\dot{\gamma_0}=\sum_{j=1}^k \dot{x}_jv_j$. Then
\[ \dot{a}_{\gamma,v_i}(r)= 2 g_{F_0}(v_i,\dot{\gamma_0}) =
2\dot{x}_i(r) \] implies $\;a_{\gamma,v_i}(r) = 2 x_i(r) + c_i\,$
and from the initial condition follows $c_i=-2x_i(0)$. Hence
$a_{\gamma,v_i}(1)=2 x_i(1)-2x_i(0) = 0$ since $\gamma$ is closed.
Then using the above formula for the parallel displacement we
obtain
\[
\tau^C_{\delta}=\tau^C_{\gamma} = \left( \begin{array}{cccc}
  1 & 0 & a_{\gamma_1}(1) & * \\[0.1cm]
  0 & Id_{ATF_0} & 0 & * \\[0.1cm]
  0 & 0 & A\circ \tau_{\gamma_1}^{F_1} \circ A^{-1} & * \\[0.1cm]
  0 & 0 & 0 & 1
\end{array}\right)
\]
This shows that the non-degenerate subspace $AT_xF_0 \subset
T_{(0,0,x)}C$ is holonomy-invariant. Hence in that case, the
cylinder is decomposable.\\[0.1cm]
Now, let $\;(F,g_F)=(F_1,g_{F_1}) \times \dots \times
(F_k,g_{F_k})\;$ be a product of irreducible Riemannian manifolds
of dimension $\geq 2$ and $x=(x_1,\dots,x_k)$. Let us denote by
\[ \fh:=\fhol_{\hat{x}}(C,g_C) \subset \fso(T_{\hat{x}}C,g_C)_{P} =
\R P\wedge AT_xF + \fso(AT_xF,g_0) + \R P\wedge Q \] the holonomy
algebra of $(C,g_C)$ in the point $\hat{x}$, where $
\fso(T_{\hat{x}}C,g_C)_{ P} $ is the Lie algebra of the stabilizer
of $P$ in $SO_0(T_{\hat{x}}C,g_C)$. Furthermore, let $\fm_i$ be
the image of the projection of $\fh$ onto $\;(\R P \wedge
AT_{x_i}F_i)\;$ and $\fh_i$ be the image of the projection of
$\fh$ onto $\fso(AT_{x_i}F_i)$. From (\ref{P11}) we know, that $\;
\fh_i =  A \circ \fhol_{x_i}(F_i,g_{F_i})\circ A^{-1}\;$. Since
$(F,g_F)$ is supposed to have parallel spinors, $(F_i,g_{F_i})$
has parallel spinors as well. Therefore the algebras $\fh_i$, being $\fsu (m), \fsp(m), \fg_2$ or $\fspin (7)$ each, have
no center. It follows from the classification of weakly
irreducible subalgebras of $\fso(1,n)_{P }$ (cf. \cite{BI} or
\cite{G1}) that there is no coupling between the $\fm_i$ and the
$\fh_i$-part, since such couplings can only appear, if one of the
Lie algebras $\fh_i$ has a center.  Hence \bea \fh =
\left(\begin{array}{cccccc}
  0 & \fm_1 & \fm_2 & \cdots & \fm_k & 0 \\
  0 & \fh_1 & 0 & \cdots & 0 & *\\
  0 & 0 & \fh_2 & \cdots & 0 & *\\
  \vdots & \vdots &\vdots
  &\vdots & \vdots &\vdots \\
  0 & 0 & 0& \cdots & \fh_k & *\\
  0 & 0 & 0&  \cdots  & 0 & 0
\end{array}\right) \subset \fso(T_{\hat{x}}C,g_C)_P \label{P16}
\eea  Now, we will show that $\fm_i \not = 0$.  Note, that by the
Ambrose-Singer-Theorem, $\fh$ is spanned by all elements of the
form
\[ (\tau^C_{\delta})^{-1} \circ R^C_{\delta(1)}(X,Y) \circ
\tau^C_{\delta} \] where $\delta: [0,1] \longrightarrow C$ runs
over all curves in $C$ starting in $\hat{x}$ and $X$ and $Y$ over
all vectors in $T_{\delta(1)}C$. Let us consider the special curve
$\delta(r)=(0,0,\gamma(r))$ where $\gamma: [0,1] \longrightarrow
F$ is a geodesic in $(F,g_F)$ starting in $x$ . Then for the
parallel displacement of a vector $Z\in T_xF$ hold by Proposition
\ref{Prop-parallel}
\[ \tau^C_{\delta}(AZ) = a_{\gamma,Z}(1) P + A(\tau_{\gamma}^FZ). \]
From Proposition \ref{Prop-curvatureFC} follows for vectors $X,Y
\in TF$ (recall $P \perp H^{-1} TF$):
\[ R^C(X,Y)P=0 \qquad \mbox{and} \qquad R^C(X,Y)H^{-1} \vert_{TF} =
H^{-1}R^F(X,Y) \vert_{TF}.\] Since $A=H^{-1}$ this implies
\[ R^C_{\delta(1)}(X,Y)\tau_{\delta}^C(AZ) = A R^F_{\delta(1)}(X,Y)
\tau^F_{\gamma}(Z).\] Let us denote by $\gamma^{-}$ the inverse
curve $\gamma^{-}(r)=\gamma(1-r)$. Then with
$\hat{Z}:=R^F_{\delta(1)}(X,Y) \tau^F_{\gamma}(Z)$ we obtain (recalling $R(X,Y)P=0$)
\[ (\tau_{\delta}^C)^{-1} \circ R^C_{\delta(1)}(X,Y)\circ
\tau_{\delta}^C(AZ) = a_{\gamma^-,\hat{Z}}P + A \tau_{\gamma^-}^F
R^F_{\gamma(1)}(X,Y)\tau_{\gamma}^F(Z) \] where
$a_{\gamma^-,\hat{Z}}$ is given by the differential equation
\[ \dot{a}_{\gamma^-,\hat{Z}}(r) = 2g_F\Big( \dot{\gamma}^-(r),
\tau^F_{\gamma^-|_{[0,r]}}R^F_{\gamma(1)}(X,Y)\tau^F_{\gamma}(Z)\Big)
\qquad \qquad (*)
\]
with the initial condition $a_{\gamma^-,\hat{Z}}(0)=0$. Since
$\gamma^-$ is a geodesic in $F$, the function $\;g_F\Big(
\dot{\gamma}^-(r),
\tau^F_{\gamma^-|_{[0,r]}}R^F_{\gamma(1)}(X,Y)\tau^F_{\gamma}(Z)\Big)\;$
is constant, hence the solution of the initial value problem (*)
is
 \bean
 a_{\gamma^-,\hat{Z}}(r) &=& 2 g_F( \dot{\gamma}^-(0),
R^F_{\gamma(1)}(X,Y)\tau^F_{\gamma}(Z)) \cdot r \\
&=& -2g_F(\dot{\gamma}(1), R^F_{\gamma(1)}(X,Y)\tau^F_{\gamma}(Z))
\cdot r \\
&=& 2
g_F(R^F_{\gamma(1)}(X,Y)\dot{\gamma}(1),\tau^F_{\gamma}(Z))\cdot r
\eean Now, assume that $\fm_i = 0$. Then
$R^{F_i}_{\gamma_i(1)}(X,Y)(\dot{\gamma_i}(1))$ vanish for all
vectors $X,Y \in T_{\gamma_i(1)}F_i$. In particular, this implies
that the sectional curvature $K^{F_i}_E(\gamma_i(1))$ of $F_i$ in
the point $\gamma_i(1)$ in direction of any 2-dimensional subspace
$E =
span(\tau^{F_i}_{\gamma_i}(v),\tau^{F_i}_{\gamma_i}(\dot{\gamma_i}(0))$
with $v\in T_{x_i}F_i$ vanishes. The same argument applies to any
point $\gamma(r)$  of $\gamma((0,1])$, hence we obtain taking the
limit $r \to 0$ that $K^{F_i}_{span(v,\dot{\gamma_i}(0))}(x_i) =
0$. If we take all geodesics $\gamma_i$ starting from $x_i$ it
follows that all sectional curvatures of $F_i$ in the point $x_i$
vanish. Now, since the holonomy groups of different points are
conjugated, with $\fm_i$ the projection of
$\fhol_{\hat{y}}(C,g_C)$ onto the part $\R P \wedge AT_{y_i}F_i
\subset \fso(T_{\hat{y}}C,g_C)_P$ for any other point $y\in F$
vanishes too. Then applying the same argument to $y$ we obtain,
that the sectional curvature of $F_i$ vanishes everywhere, hence
$F_i$ has to be flat, which is a contradiction since
$(F_i,g_{F_i})$ is irreducible and of dimension $\geq 2$. Hence we
have proven that the projection $\fm_i \not = 0$  for all
components $F_i$ of $F$. \\
It remains to prove that $\fm_i=AT_{x_i}F_i$ for all $i$. Formula
(\ref{P16}) shows that $ \, [\fh_i, P\wedge \fm_i] \subset P
\wedge \fh_i(\fm_i)\; $ hence, $\fh_i(\fm_i) \subset \fm_i$. But
by assumption, $\fh_i$ acts irreducible on $AT_{x_i}F_i$, hence
$\fm_i = AT_{x_i}F_i$.  This completes the proof of the second
statement. \hfill \qed

\begin{Corollary}
Let an indecomposable, non-irreducible Lorentzian holonomy representation
$\mathfrak{R}$ with a fixed spinor be given.
Then there is a bbc manifold with holonomy representation $\mathfrak{R}$.
\end{Corollary}

\v {\bf Proof.} First pick a complete manifold $(F, g_F)$ having the screen bundle holonomy (cf. \cite{tl}) of $\mathfrak{R}$ as its holonomy representation. Then from any (possibly trivial) Codazzi tensor on $(F, g_F)$ with norm uniformally bounded from below one can construct a Codazzi tensor $H$ on the warped product over $(F,g_F)$ with positive eigenvalues uniformally bounded away from zero. Now $C[F;H]$ is a bbc manifold with holonomy representation $\mathfrak{R}$.  \hfill \qed

\bigskip

\ \\
\section{Examples}\label{sec-Examples}

Our construction of Riemannian manifolds with Codazzi spinors and
of Lorentzian manifolds with special holonomy (cf. Theorem
\ref{Th-Codazzi-spinor1} and Theorem \ref{Theorem-hol1}) is based
on the existence of Codazzi tensors on Riemannian manifolds with
parallel spinors. Let us discuss some examples for that.\\

{\bf Example 1}\\
On the flat space $\R^n$ the endomorphism
\[ T^{\R^k}_h(X) = \nabla^{\R^k}_X (grad(h)) = X(\partial_1h, \dots, \partial_kh))\]
is a Codazzi tensor for any function $h$ on $\R^k$ and every
Codazzi tensor is of this form (cf. \cite{Ferus}). Proposition
\ref{Prop-curvatureFC} shows, that the cylinder $C[F;H]$ is flat
for any Codazzi tensor $H$ on the warped product that
is constructed out of $T$ (cf. Corollary \ref{Cor-Codazzi2}).\\

{\bf Example 2}\\
Let $(F_1,g_{F_1})$ be a complete simply connected irreducible
Riemannian spin manifold with parallel spinors and $(F,g_F)$ its
Riemannian product with a flat $\R^k$. Then $(F,g_F)$ is complete
and has parallel spinors. Let $H$ be a Codazzi tensor on the
warped product $\R \times_{e^{-2s}} F$ constructed out of the
Codazzi tensor $\lambda Id_{F_1} + T_h^{\R^k}$ of $F$, where
$T_h^{\R^k}$ is taken from example 1. Then the cylinder $C(F;H)$
is globally hyperbolic and decomposable with special holonomy
\[ Hol(F_1,g_{F_1})\;\ltimes\; \R^{dim F_1}.\]  \\

{\bf Example 3}\\
Let us consider the metric cone
\[ (F^{n-1},g_F) := (\R^+ \times N, dr^2 + r^2 g_N) \]
where $(N,g_N)$ is simply connected and a Riemannian
Einstein-Sasaki manifold, a nearly K\"ahler manifold, a 3-Sasakian
manifold or a 7-dimensional manifold with vector product. Then
$(F,g_F)$ is irreducible and has parallel spinors (but fails to be
complete). Furthermore,  $T:= \nabla^F
\partial_r$ is a Codazzi tensor on $(F,g_F)$, since $\partial_r$ lies in the kernel
of the curvature endomorphism. Theorem \ref{Theorem-hol1} shows,
that the cylinder $C[F;H]$, where $H$ is constructed out of $T$,
has special holonomy
\[ Hol(C,g_C) \simeq G \ltimes \R^{n-1} \]
where
\[ G= \left\{ \begin{array}{lll}
 SU((n-1)/2) & \qquad &  \mbox{if $\quad N$ is Einstein-Sasaki} \\
 Sp((n-1)/4) & \qquad & \mbox{if $\quad N$ is 3-Sasakian}\\
 G_2         & \qquad & \mbox{if $\quad N$ is nearly K\"ahler}\\
 Spin(7)     & \qquad & \mbox{if $\quad N$ 7-dimensional with vector product}
 \end{array}\right.
 \]
\ \\

{\bf Example 4}\\
Let $(F,g_{F})=(F_1,g_{F_1}) \times \dots \times \dots
(F_k,g_{F_k})$ be a Riemannian product of simply connected
complete irreducible Riemannian manifolds with parallel spinors.
Let $T$ be the Codazzi tensor $T= \lambda_1 \1_{F_1} + \dots
\lambda_k \1_{F_k}$ and $H$ constructed out of $T$ as mentioned in
Corollary \ref{Cor-Codazzi2}. Then $C[F;H]$ is globally
hyperbolic, weakly irreducible and the holonomy group is
isomorphic to
\[ \left(Hol(F_1,g_{F_1}) \times \dots \times
Hol(F_k,g_{F_k})\right) \; \ltimes \; \R^{dim F}. \]

\ \\
{\bf Example 5.  Eguchi-Hansen space}\\
We will show that there is no nontrivial Codazzi tensor on the
Eguchi-Hansen space, which is an example of a complete,
irreducible Riemannian 4-manifold with holonomy $SU(2)$, hence
with 2 linearly independent parallel spinors.\\
The Eguchi-Hansen space $\overline{EH}$ is the manifold $T \S^2$
equipped with a metric which is ALE (asymptotically locally
Euclidean). More exactly, there is a compact set whose complement
converges with the metric distance to $\R^4 / \Z_2$ in a certain
way (cf. \cite{jo1}). The complement of the zero section of
$T\S^2$ is isometric to
 $((a, \infty) \times \S^3)/\{ \pm 1 \}$ with the metric

$$h_a = (1-(\frac{a}{r})^4)^{-1} dr^2 + r^2 ((\sigma^x)^2 + (\sigma^y)^2) +
r^2(1-(\frac{a}{r})^4) (\sigma ^z)^2 ,$$

\v where $r$ is the parameter of $(a, \infty) $ and $(\sigma^x,
\sigma^y, \sigma^z)$ is the standard basis of left-invariant
1-forms on $\S^3$, i.e. $d \sigma^x = 2 \sigma^y \wedge \sigma^z$
(and cyclic permutations). (Note that the metric is in fact
invariant under the multiplication by $-1$ in $\R^4$!). We denote
this complement of the zero section by $EH$ as, vice versa, $
\overline{EH} $ is the completion of $EH$ by gluing in an $\S^2$
at $r \rightarrow a$ (cf. \cite{ra}, \cite{eh}). One can extend
$h_a$ to a complete metric on $\overline{EH}$ if $a>0$.

\begin{Theorem}\label{Th-Codazzi-EH}
Let $W$ be a Codazzi tensor on $EH$. Then $W$ is a constant
multiple of the identity. Moreover, on the warped products $\R \times_{e^{ks}} EH$, all Codazzi tensors are simple.
\end{Theorem}

{\bf Proof.} Let $(\sigma_x, \sigma_y, \sigma_z)$ be the basis of
left-invariant vector fields on $\S^3$ dual to $(\sigma^x,
\sigma^y, \sigma^z)$. Then  $[\sigma_x, \sigma_y] = - 2 \sigma_z$
(and cyclic permutations).  Written in the standard spherical
basis $(r
\partial_r, \sigma_x, \sigma_y, \sigma_z)$, the
metric becomes
$$ h_a =r^2
\left(
\begin{array}{cccc}
(1-(\frac{a}{r})^4)^{-1}& 0& 0& 0 \\
0& 1& 0& 0 \\
0& 0& 1& 0 \\
0& 0& 0& (1- (\frac{a}{r})^4)
\end{array} \right)
$$
The transformation matrix between $(r
\partial_r, \sigma_x, \sigma_y, \sigma_z)$ and the canonical basis
of the Euklidian coordinates on $\C^2 = \R^4$ at a point $(x_1,
x_2, x_3, x_4)$ is
$$ M =
\left(
\begin{array}{cccc}
x_1& x_2& x_3& x_4 \\
x_2& -x_1& x_4& -x_3 \\
x_3& -x_4& -x_1& x_2 \\
x_4& x_3& -x_2& -x_1
\end{array} \right)
$$
For short, we denote by $f$ the function
\[ f(r) = \Big(1-\Big(\frac{a}{r}\Big)^4 \Big)^{1/2}. \]
Note, that for this function
\[ f' + r^{-1}f = 2r^{-1}f^{-1}-r^{-1}f =: \gamma(r) .
\]
Then
\[ \sigma^0 = f^{-1}dr \;, \;\; \sigma^1 = r \sigma^x \; ,\; \; \sigma^2 = r
\sigma^y  \;, \;\;\sigma^3 = rf\sigma^z \] is an orthogonal
coframe for the metric $h_a$ and
\[ e_0=f \partial_r \; , \;\; e_1= r^{-1}\sigma_x \;,\; \; e_2 =
r^{-1}\sigma_y \;,\; \; e_3 = r^{-1}f^{-1}\sigma_z \] is the dual
orthonormal basis for $h_a$. For the commutators we obtain \bean
\, [e_0,e_1]& =& -
r^{-1}f \,e_1 \\ \,[e_0,e_2] &=& - r^{-1}f \,e_2 \\
\, [e_0,e_3] &=& - (f'+ r^{-1}f) \,e_3 \; = \; -\gamma\,e_3 \\
\, [e_1,e_2]& = & -2r^{-1}f \,e_3 \\
\, [e_1,e_3] &=&
\, + 2r^{-1}f^{-1}\, e_2 \\
\, [e_2,e_3] &=&  -2r^{-1}f^{-1}\, e_1 \eean
Then using the Koszul formula
\[ 2 \langle \nabla_{e_i}e_j,e_k \rangle = -\langle e_i,[e_j,e_k]\rangle
- \langle e_j, [e_i,e_k] \rangle + \langle e_k,[e_i,e_j] \rangle
\] we obtain for the Levi-Civita connection of $h_a=\langle
\cdot,\cdot \rangle$
\bean
\nabla_{e_0}e_k &=& 0 \;\;\;\quad k=0,1,2,3\\[0.2cm]
\nabla_{e_1} e_0 &= &+ r^{-1}f e_1  \qquad \qquad \qquad \nabla_{e_2} e_0 \; = \; + r^{-1} f e_2\\
\nabla_{e_1} e_1 &=& - r^{-1} f e_0 \qquad \qquad \qquad \nabla_{e_2} e_1 \; = \; + r^{-1} f e_3\\
\nabla_{e_1} e_2 & = & - r^{-1} f e_3 \qquad \qquad \qquad \nabla_{e_2} e_2 \; =\; - r^{-1} f e_0\\
\nabla_{e_1} e_3 &=& + r^{-1} f e_2 \qquad \qquad \qquad  \nabla_{e_2} e_3 \; =\;  - r^{-1} f e_1 \\[0.2cm]
\nabla_{e_3} e_0 &=& + ( r^{-1} f  + f') \;e_3 \;= \; +\gamma\;e_3\\
\nabla_{e_3} e_1 &=& + ( r^{-1} f  - 2r^{-1} f^{-1})\; e_2 \; = -\gamma\;e_2\\
\nabla_{e_3} e_2& =& - ( r^{-1} f  - 2 r^{-1} f^{-1})\; e_1 \; = +\gamma \; e_1\\
\nabla_{e_3} e_3 &=& - ( r^{-1} f  + f')\; e_0 \; = \; -\gamma \;
e_0 \eean
We consider $W$ in the basis $e_a$ mentionned above. The first
observation is that w.r.o.g the entries of the associated matrix
only depend on $r$, not on the spherical variables. That can be
seen as follows: Assume that there is a Codazzi tensor field $W$
whose entries might depend on all variables. Then consider
$$\tilde{W} := \int_{\S^3} (\1 \times L_q)^* W dq $$
where $L_q$ is the left action of the sphere on itself. Remind
that the Berger metrics are left-invariant, thus the
diffeomorphisms $\1 \times L_q$ are isometries of $EH$. Therefore
the endomorphisms $(\1 \times L_q)^* W$ are Codazzi. As the
Codazzi equation is linear in the tensor field, it commutes with
the integral, and $\tilde{W}$ is Codazzi. As the frame $e_a$ is
left-invariant, the entries of $\tilde{W}$ cannot depend on the
spherical coordinates any more but only on $r$. Up to know we know
that w.r.o.g. we have a tensor field of the form
$$ W :=
\left(
\begin{array}{cccc}
A& B& C& D \\
B& E& F& G \\
C& F& H& I \\
D& G& I& J
\end{array} \right)
$$
with $A,B,C,D,E,F,G,H,I,J$ real functions depending on $r$. The
second observation is that, for $U$ being the endomorphism which
exchanges $e_1$ and $e_2$ and is the identity on the orthogonal
complement of $\{ e_1, e_2 \}$, if $W$ is a Codazzi tensor, the
endomorphism $U \o W \o U$ is a Codazzi tensor as well. Thus the
endomorphism fields $W \pm U \o W \o U$ are Codazzi tensors as
well which are symmetric resp. antisymmetric under the conjugation
by $U$. We will show that the $U$-antisymmetric part vanishes
necessarily while the $U$-symmetric part has to be a constant
multiple of the identity. \\
The $U$-antisymmetric part has
$A=D=F=J=0$ and $B=-C$, $E=-H$, $G= -I$. Therefore it looks like
$$ W :=
\left(
\begin{array}{cccc}
0& B& -B& 0 \\
B& E& 0& G \\
-B& 0& -E& -G \\
0& G& -G& 0
\end{array} \right)
$$
Now we consider the Codazzi equation for the vectors $e_1, e_2$:
\bean
0 &=& \nabla_{e_1} (W(e_2)) - \nabla_{e_2} (W(e_1)) - W ([e_1, e_2])\\
  &=& \nabla_{e_1} (-Be_0 -E e_2 - G e_3) - \nabla_{e_2} (B e_0 + E e_1 + G e_3 ) - W ( -2 \frac{f}{r} e_3)\\
  &=& \frac{f}{r} ( (-B e_1  +E e_3 - G e_2) - (B e_2 + E e_3 - G e_1) + 2 ( G e_1 - G e_2 ))\\
  &=& \frac{f}{r} \big( (-B + 3G) e_1 + (-B - 3G) e_2   \big)
\eean
which implies $ B= G=0$. Thus $W$ looks like
$$ W :=
\left(
\begin{array}{cccc}
0& 0& 0& 0 \\
0& E& 0& 0 \\
0& 0& -E& 0 \\
0& 0& 0& 0
\end{array} \right)
$$
The same procedure applied to the pair $e_1, e_3$ leads to
\bean
0 &=& \nabla_{e_1} (W(e_3)) - \nabla_{e_3} (W(e_1)) - W ([e_1, e_3])\\
  &=& -E \nabla_{e_3} e_1 - W (2 r^{-1}f^{-1} e_2) \\
  &=& E \g e_2 + E 2 r^{-1} f^{-1} e_2
  \eean
which implies $E=0$ as $\g = 2r^{-1}f^{-1}-r^{-1}f$.\\
On the other hand, the $U$-symmetric part has $C=B, H=E, G=I$,
i.e. that
$$ W :=
\left(
\begin{array}{cccc}
A& B& B& D \\
B& E& F& G \\
B& F& E& G \\
D& G& G& J
\end{array} \right)
$$
Now consider the Codazzi equation for the vectors $e_1, e_2$:
\bean
0 &=& \nabla_{e_1} (W(e_2)) - \nabla_{e_2} (W(e_1)) - W ([e_1, e_2])\\
  &=& \nabla_{e_1} (B e_0 + F e_1 + K e_2 + G e_3) - \nabla_{e_2} (B e_0 +
  K e_1 + F e_2 + G e_3 ) - W ( -2 \frac{f}{r} e_3)\\
  &=& \frac{f}{r} ( (B e_1 - F e_0 - K e_3 + G e_2) - (B e_2 + K e_3 - F e_0 - G e_1)
   + 2 (D e_0 + G e_1 + G e_2 + Je_3))\\
  &=& 2D e_0 + (B + 3G) e_1 + (-B + 3G) e_2  + (2J-2K) e_3
\eean
which implies $D= B= G=0$. Therefore $W$ has necessarily the form
$$ W :=
\left(
\begin{array}{cccc}
A& 0& 0& 0 \\
0& K& F& 0 \\
0& F& K& 0 \\
0& 0& 0& K
\end{array} \right)
$$
But then the Codazzi equation applied to the vectors $e_1, e_3$
reads
\bean
0 &=& \nabla_{e_1} (W(e_3)) - \nabla_{e_3} (W(e_1)) - W ([e_1, e_3])\\
&=& K \nabla_{e_1} e_3 - K \nabla_{e_3}e_1 - F \nabla_{e_3} e_2 - 2 r^{-1}f^{-1} (Fe_1 + Ke_2)\\
&=& F( -\g - 2 r^{-1}f^{-1}) e_1 \eean
which can never be satisfied unless $F=0$ as $\g =
2r^{-1}f^{-1}-r^{-1}f$. Thus $W$ must have the form
$$ W :=
\left(
\begin{array}{cccc}
A& 0& 0& 0 \\
0& K& 0& 0 \\
0& 0& K& 0 \\
0& 0& 0& K
\end{array} \right)
$$
Finally, we consider the vectors $e_0, e_1$ and the vectors $e_0,
e_3$. The Codazzi equation of the first pair is equivalent to the
condition $e_0 (K) = \frac{f}{r} (A-K)$, the Codazzi equation of
the second pair is equivalent to $e_0(K) = \g (A-K)$. Combined
this yields $(\frac{f}{r} - \g )(A-K) =0$, thus $A=K$ and $e_0 (K)
=0$, thus we have shown that the $U$-symmetric part is a constant
multiple of the identity on $EH$ which completes the proof of the
first part of the Theorem. The second part is proven analogously: First we assume the existence of a nonzero homothetic vector field on $EH$, i.e. of a vector field $V$ with $\nabla_{\cdot} V = c \cdot \1$. By integrating over $\S^3$ we show that all coefficients of $V$ w.r.t. the left-invariant basis above cannot depend on the spherical variables but only on $r$. Thus $V = Ae_0 + B e_1 + C e_2 + D e_3$ with $A,B,C,D$ real functions depending only on $r$. Then we compute

$$c e_1 = \nabla_{e_1} V = A \nabla_{e_1} e_0 + B \nabla_{e_1} e_1 + C \nabla_{e_1} e_2 + D \nabla_{e_1} e_3$$

\v and as only $\nabla_{e_1} e_0 $ is collinear to $e_1$, we conclude that $B=C=D=0$ and $c = \frac{f}{r}A$. But the same procedure for $e_3$ instead of $e_1$ gives $c= \g A$. Thus both equations together can be satisfied only for $A=0$ as $\g =2r^{-1}f^{-1}-r^{-1}f$. Thus there is no nonzero vector field $V$ on $EH$ with $\nabla_{\cdot} V = c \cdot \1$. Therefore Corollary \ref{Cor-Codazzi1} implies that on $\R \times_{e^{ks}} EH$ there are no non-simple Codazzi tensors.  \hfill \qed

\section{Acknowledgements}

\v Both authors wish to express their gratitude for the hospitality of the Erwin-Schr\"odinger-Institut in Vienna where substantial parts of this article were obtained.


\small{

}


\begin{thebibliography}{99}

\bibitem{bgm}
Ch.~B\"ar, P.~Gauduchon, A.~Moroianu: Generalized cylinders in
semi-Riemannian and spin geometry. Math. Zeitschrift 249 (2005),
545--580.

\bibitem{baum1}
H.~Baum: Spin-Strukturen und Dirac-Operatoren \"uber
pseudo-Riemannschen Mannigfaltigkeiten. Teubner-Texte zur
Mathematik. Bd. 41, Teubner-Verlag Leipzig, 1981.

\bibitem{h1}
H.~Baum: Riemannian manifolds with imaginary Killing spinors,
Annals of Global Analysis and Geometry vol.7, no 2 (1989),
p.141-154

\bibitem{h2}
H.~Baum: Complete Riemannian manifolds with imaginary Killing
spinors, Annals of Global Analysis and Geometry vol.7, no 3
(1989), p. 205-226

\bibitem{BFGK}
H.~Baum, T.~Friedrich, R.~Grunewald, I.~Kath: Twistors and Killing
spinors on Riemannian manifolds. Teubner-Texte zur Matheamtik, Bd.
124. Teubner-Verlag Leipzig/Stuttgart, 1991.

\bibitem{bee}
J.~K.~Beem, P.~E.~Ehrlich, K.~L.~Easley: Global Lorentzian
Geometry, 2nd edition, Marcel Decker Inc. 1996.

\bibitem{BI}
L.~Berard-Bergery, A.~Ikemakhen: On the holonomy of Lorentzian
manifolds. In: {\em Differential Geometry: Geometry in
Mathematical Physics and Related Topics}. Volume 54 of Proc.
Sympos. Pure Math., 1993, p 27--40.

\bibitem{dc}
A.~Derdzinski, Ch.-L. Sen: Codazzi tensor fields, curvature and
Pontryagin forms, Proc. Lond. Math. Soc. (3), 47 (1983), p.15-26.

\bibitem{eh}
T.~Eguchi, A.~L.~Hanson: Asymptotically flat self-dual solutions
to Euclidean gravity. Physics Letters 74B, (1978), pp. 249-251.

\bibitem{Ferus}
D.~Ferus: A remark on Codazzi tensors in constant curvature
spaces. Lecture Notes in Math. 838, Springer, 1981, p. 257.

\bibitem{G1}
A.~Galaev: The space of curvature tensors for holonomy algebras of
Lorentzian manifolds. Diff. Geom. and its Appl. 22 (2005), 1-15.

\bibitem{G2}
A.~Galaev: Metrics that realize all types of Lorentzian holonomy
algebras. arXiv:mathDG/ 0502575, 2005.

\bibitem{h}
N.~Hitchin: Harmonic spinors, Advances in Mathematics 14 (1974),
p.1-55.

\bibitem{jo1}
D.~Joyce: Asymptotically Locally Euclidean metrics with holonomy
$SU(m)$. arXiv: math.AG/ 9905041, 1999.

\bibitem{jo2}
D.~Joyce: Compact manifolds with special holonomy. Oxford Science
Publications, 2000.

\bibitem{kn}
S.~Kobayashi, K. Nomizu: Foundations of differential geometry,
vol.1, Interscience Publishers, 1969.

\bibitem{tl}
Th.~Leistner: Holonomy and Parallel Spinors in Lorentzian
Geometry, PhD thesis, Humboldt University Berlin 2003. Logos
Verlag 2004.

\bibitem{tl1}
Th.~Leistner: Torwards a classification of Lorentzian holonomy
groups. arXiv:mathDG/ 0305139, 2003.

\bibitem{tl2}
Th.~Leistner: Torwards a classification of Lorentzian holonomy
groups. Part II. arXiv: mathDG/0309274, 2003.

\bibitem{m}
O.~M\"uller: The Cauchy Problem of Lorentzian Minimal Surfaces in Globally Hyperbolic Manifolds, accepted for publication by: Annals of Global Analysis and Geometry, 2005

\bibitem{ONeill}
B. O'Neill: Semi-Riemannian Geometry. Academic Press, Inc. 1983.

\bibitem{ra}
P.~Ramacher: Geometric and analytic properties of families of
hypersurfaces in Eguchi-Hanson space. J. Geom. Phys. 44 (2003);
407--474.


\end{thebibliography}
\end{document}